%
\magnification = 1200
\raggedright

\raggedbottom  
\input amssym.def
\input amssym.tex


\font\bigtenrm = cmr10 scaled \magstep1

\def\bwedge{\bigwedge}
\def\cwedge{\bigwedge}

\def\symmdif{\bigtriangleup}
\def\tensor{\otimes}
\def\Club{\rm{Club}}

\def\Sk{\rm{Sk}}

\def\order{\prec}

\def\gl{\rm{gl}}

\def\Sk{\rm{Sk}}

\def\bd{\rm{bd}}                      
\def\cf{\rm{cf}}
\def\into{\rightarrow}

\def\rest{\upharpoonright}  
\def\eop{$\bigstar$}  
\def\1st{1$^{\hbox{st}}$}
\def\2nd{2$^{\hbox{nd}}$}
\def\card#1{\vert{#1}\vert} 

\def\cf{\mathop{\rm cf}}
\def\ref{\mathop{\rm ref}}  
\def\nacc{\mathop{\rm nacc}}
\def\acc{\mathop{\rm acc}}
\def\implies{\Longrightarrow}
\def\iff{\Leftrightarrow}
\def\satisfies{\vDash}

\def\deq{\buildrel{\rm def}\over =}
\def\otp{\mathop{\rm otp}}

\def\to{\longrightarrow}
\def\mod{\mathop{\rm mod}}
\def\elementary{\order}

\def\CC{{\cal C}}
\def\DD{{\cal D}}

\def\PP{{\cal P}}
\def\RR{{\cal R}}


\def\today{\ifcase\month\or
   January\or February \or March \or April \or May \or June \or
   July \or August \or September \or October \or November \or December \fi
   \space \number\day, \number\year}


\centerline{{\bigtenrm On squares,
outside guessing of clubs and
 $I_{<f}[\lambda]$}
}

\medskip

\centerline{Mirna D\v zamonja and Saharon Shelah 
\footnote{${}^1$}{\sevenrm Authors partially
supported by the Basic Research Foundation Grant number 0327398
administered by the Israel Academy of Sciences and Humanities.
The first author thanks the Hebrew University and the Lady Davis Foundation
for the Forchheimer Postdoctoral Fellowship.
For easier future
reference, note that this is publication [DjSh 562] in Shelah's bibliography.
The results presented were obtained in the period April to August 1994.
The appendix was added in December 1994.
We wish to thank Moti Gitik and Ofer Shafir for their interest and helpful comments,
as well as to James Cummings for pointing out a difficulty.}}

\centerline{\sevenrm Hebrew University of Jerusalem}

\centerline{\sevenrm Institute of Mathematics}

\centerline{\sevenrm  91904 Givat Ram}

\centerline{\sevenrm  Israel}

\centerline{\sevenrm dzamonja@math.huji.ac.il  \&  shelah@math.huji.ac.il}



\centerline{to be published in Fundamenta Mathematicae 148(3)}

\vskip 2 true cm

\centerline{\bf Abstract}

\smallskip

\smallskip

\smallskip

Suppose that $\lambda=\mu^+$ and $\mu$ is singular.
We consider two aspects of the square property on subsets of $\lambda$.
First, we have results which 
show e.g. that for $\aleph_0<\kappa=\cf (\kappa)<\mu$,
the equality $\cf([\mu]^{\le\kappa},\subseteq)=\mu$
is a sufficient condition 
for the set of elements of $\lambda$
whose cofinality is $\le\kappa$, to
be split into the union of $\mu$
sets with squares. 
Secondly, we
introduce a certain weak version of the square property and prove that
if $\mu$ is a strong limit, then
this weak square property holds
on $\lambda$ without any additional assumptions.

In the second section we start with
two universes $V_1\subseteq V_2$ of set theory,
and a regular cardinal $\kappa$ in $V_1$ 
such that
the cofinality of $\kappa$
in $V_2$ is $\theta<\kappa$.
Assume $\kappa^+$ is preserved and
$\kappa$ is inaccessible in $V_1$ with $2^\kappa=\kappa^+$.
We show that then
there is an unbounded subset $C$ of $\kappa$ in $V_2$, such that
for every club $E$ of $\kappa$ in $V_1$, the difference $C\setminus E$ is
bounded. 
We have further results of a similar flavor.
Some of our results were independently obtained by Moti Gitik, using different
methods.

In the third section we consider the connection between the ideal 
$I[\lambda]$ and
the notions of square and weak square. We show that these notions are a
part of a larger family of properties which all can be introduced through
a single definition of $I_{<f}[\lambda]$ by changing the parameter
$f$. We discuss further properties of $I_{<f}[\lambda]$ and some other
similarly defined notions. We have further results on $I[\lambda]$
in the last section. 

\vfill

\eject
\baselineskip=24pt

{\bf \S0. Introduction.} The problems studied in this paper come naturally
in the study of cardinal arithmetic. The notions involved, like
the ideal $I[\lambda]$, decomposition into sets with squares and club guessing
have been extensively investigated and applied by the second author in
[Sh -g] and related papers, both before and after [Sh -g]. 

In [Sh 351,\S4]
and [Sh 365, 2.14] it is shown that if $\mu$ is a regular cardinal, then
$\{\alpha<\mu^+:\,\cf(\alpha)<\mu\}$ can be written as the union of
$\mu$ sets on which there are squares.
In 1.1 of this paper it is shown that
for a singular cardinal $\mu$ and $\aleph_0<\cf (\kappa)
=\kappa<\mu$,
if $\cf\bigl([\mu]^{\le\kappa},\subseteq\bigr)=\mu$,
then $\{\alpha<\mu^+:\,\cf(\alpha)\le\kappa\}$ is the union of $\mu$
sets with squares. The proof is an application of [Sh 580].
The
present result improves [Sh 237e, 2] for
a singular $\mu$, as [Sh 237e, 2] which had the same conclusion
and assumed $\mu^{\le\kappa}=\mu$. It also implies that under the
assumptions of 1.1, the set $\{\alpha<\mu^+:\,\cf(\alpha)\le\kappa\}$ is an
element of $I[\mu^+]$--a fact which also follows from [Sh 420, 2.8.].
Here $I[\mu^+]$ is the ideal introduced in [Sh 108] or [Sh 88a].

Also in the first section is a theorem which shows that if $\mu$ is a
singular strong limit, then there is a weak version of the square
principle, which we call square pretender, such that ``many''
elements of $\mu^+$ have a club on which there is a square pretender.
Moreover, all square pretenders in question can be enumerated in
type $\mu$.

Suppose $\kappa$ is an inaccessible such that
$2^\kappa=\kappa^+$ and we change its
cofinality to $\theta<\kappa$, so that $\kappa^+$ is preserved.
Then there is an unbounded subset $C$ of $\kappa$ in the extension, such that
for every club $E$ of $\kappa$ in the ground model, $C\setminus E$
is bounded. This is one of the results of \S2. We have further results of this
nature, and with different assumptions.
We shall refer to this type of results as to
``outside guessing of clubs". 
Results on guessing clubs are reasonably well known (see [Sh -g], [Sh -e]).
When Moti Gitik told the second author about his result quoted in A below,
the second author was reminded of his earlier result quoted in B below,
which was done in the preprint [Sh -e], for a given club guessing.
Note the connection between A and B via generic ultrapowers. The results of the form A are wider, as they 
also apply to presaturated ideals. It was then natural
to try to prove such results using club guessing, and this is exactly what
is done here.
We quote the theorems we referred to as A and B above:

\medskip

{\bf Theorem A (Gitik).}{\it [Gi1, 2.1.].\/}
Let $V_1\subseteq V_2$ be two models of $ZFC$.
Let $\kappa$ be a regular cardinal of $V_1$ which changes its cofinality to
$\theta$ in $V_2$. Suppose that in $V_1$ there is an almost increasing
(mod nonstationary) sequence of clubs of $\kappa$ of length $\chi$,
with $\kappa^+\le\chi$ such that every club of $\kappa$ of $V_1$
is almost included in one of the clubs of the sequence. Assume that
$V_2$ satisfies the following:

\item{(1)} $\cf (\chi)\ge (2^\theta)^+$ or $\cf (\chi)=\theta$.
\item{(2)} $\kappa\ge (2^\theta)^+.$

Then in $V_2$ there exists a sequence
$\langle\tau_i:\,i<\theta\rangle$
cofinal in $\kappa$, consisting of ordinals of cofinality
$\ge \theta^+$ so that every club of $\kappa$ of $V_1$
contains a final segment of $\langle \tau_i:\,i<\theta\rangle.$

\medskip

{\bf Theorem B (Shelah).}{\it [Sh -e III 6.2.B old version]
=[Sh -e IV 3.5 new version].\/} Let $\lambda$ be regular $>2^\kappa$
and $\kappa$ regular uncountable.
Suppose that $S\subseteq\{\delta<\lambda:\,\cf (\delta)=\kappa\}$
is stationary and $I$ is a normal ideal on $\lambda$ such that $S\notin I$.
If $I$ is $\lambda^+$-saturated, {\it then\/} we can find a sequence
(called a club system)
$\langle C_\delta:\,\delta\in S\rangle$ such that each $C_\delta$
is a club of $\delta$ of order type $\cf (\delta)$, and
for every club $C$ of $\lambda$ the set
$\{\delta\in S:\,C_\delta\setminus C\hbox{ is unbounded in }\delta\}\in I$.

\medskip

The proof of this theorem in fact gives that for every $S$ stationary in
$\lambda$:

$(\ast)_{\lambda,S}$
There exists $ S_1\subseteq S$ stationary such that we can find
a club system $\langle C_\delta:\,\delta\in S_1\rangle$ such that
$$(\forall C\hbox{ a club of } \lambda)
(\{\delta\in S_1:\,\delta>\sup (C_\delta\setminus C)\}\hbox{ is
not stationary.})$$

In the third section of the paper we unify the notions of square, weak square,
silly square and $I[\lambda]$ by a single definition of $I_{<f}[\lambda]$,
where $f$ is a parameter. We consider various properties of $I_{<f}[\lambda]$.

The last section of the paper is an appendix added after the paper
was submitted. In it we prove two further theorems on $I[\lambda]$.

Before going on to the first section, we shall review some notation
and conventions commonly used in the paper.

\medskip

{\bf Notation 0.0}(0)
Suppose that $\gamma\ge\theta$ and $\theta$ is a regular
cardinal.
Then 
$$S_{<\theta}^\gamma=\{\delta <\gamma:\aleph_0\le\cf(\delta) <\theta\}.
$$

More generally, we use
$S^\gamma_{{\rm r}\, \theta}$ for ${\rm r}\in
\{<,\le,=,\neq,>,\ge\}$ to describe
$$S^\gamma_{{\rm r}\,\theta}=\{\delta<\gamma:\,\aleph_0\le\cf(\delta)\,\&\,
\cf(\delta)\,{\rm r}\,\theta\}.$$

We use $S^\gamma_1$ to denote the set of successor ordinals below $\gamma$.

(1) For us an ``inaccessible'' is simply a regular limit cardinal $>\aleph_0$.
Similarly to (0), we define
$$S^{\rm in}_\lambda=
\{\mu<\lambda:\,\mu\hbox{ is inaccessible and }<\lambda\}.$$

(2) $SING$ denotes the class of
singular ordinals, that is, all ordinals $\delta$ with $\cf(\delta)<\delta$.
$REG$ is the class of regular cardinals.

(3) For $\lambda$ a regular cardinal $>\aleph_0$,
we denote by $\Club(\lambda)$
the club filter on $\lambda$.

The ideal of non-stationary subsets of $\lambda$ is denoted by
$NS[\lambda]$.

Sometimes we also speak of the club subsets of a $\lambda$
which does not obey the above restriction, but we shall
point this out in each particular case.

(4) If $C\subseteq\lambda$, then
$$\acc(C)=\{\alpha\in C:\,\alpha=\sup(C\cap\alpha)\}\quad\hbox{and}
\quad\nacc(C)=C\setminus\acc(C).$$

(5) If $\frak A$ is a model on $\lambda$ and $a\subseteq\lambda$, then
$\Sk_{\frak A}$$(a)$ stands for the Skolem hull of $a$ in $\frak A$.

(6) In the notation $\langle H(\chi),\in, {\order}^\ast\rangle$, the
symbol ${\order}^\ast$ stands for the well ordering of $H(\chi)$.

(7) $J^{\bd}_\kappa$ is the ideal of bounded subsets of $\kappa$,
where $\kappa$ is a cardinal.
\bigskip

{\bf \S1. On the square property.}
Our first concern is an instance of decomposing
$S^{\mu^+}_{\le\kappa}\cup S^{\mu^+}_1$
into $\mu$ sets with squares, to be made more precise in a moment. We
recall the
definition of a square sequence on a set of ordinals:

\medskip

{\bf Definition 1.0.} Suppose that $S$ is a set of ordinals
and $\kappa$ is an ordinal.
The sequence $\bar C=\langle C_\delta:\,\delta\in S \rangle$
is a {\it square on\/} $S$
{\it type-bounded by\/} $\kappa$
iff the following holds for $\delta\in S$:

\item{$(a)$} $C_\delta\subseteq\delta$ is closed.

\item{$(b)$} If $\delta$ is a limit ordinal, $C_\delta$ is unbounded in
$\delta$.

\item{$(c)$} $\beta\in C_\delta\implies\beta\in S$.

\item{$(d)$}
$\beta\in {\rm acc}
(C_\delta)\implies  C_\beta=C_\delta
\cap\beta$.

\item{$(e)$} $\delta\in S\implies \otp (C_\delta)<\kappa$.


\medskip

{\bf Theorem 1.1.} Suppose that $\mu$ is singular, $\lambda=
\mu^+$, and $\aleph_0<\cf (\kappa) =\kappa<\mu$ is such that
$\cf\bigl([\mu]^{\le\kappa},\subseteq\bigr)=\mu$.

{\it Then\/} the set
$$\{\alpha<\lambda:\,\cf(\alpha)\le\kappa\}$$
is the union of $\mu$ sets with squares
which are all type-bounded by $\kappa^+$.

{\bf Proof.} It suffices to decompose $\lambda\setminus\mu$ into $\mu$ sets with squares.
We shall fix a model $\frak{A}=\langle H(\chi),\in,\elementary^{\ast}\rangle$ for some large enough
$\chi$.

For a moment, let us also fix an $a\in [\mu]^{\le\kappa}$.
We define 
$$X_a\deq\{\gamma\in[\mu,\lambda):\,{\rm cf}(\gamma)\le\kappa\,\,\&\,\,
{\rm Sk}_{\frak A}(a\cup\{\mu,\gamma\})\cap\mu=a\}.$$
We also define
$$Y_{a,\gamma}\deq {\rm Sk}_{\frak A}(a\cup\{\mu,\gamma\})\cap\gamma.$$
It can be seen that
the sets $\langle Y_{a,\gamma}\setminus\mu:\,\gamma\in X_a\rangle$ are quite close to a square sequence on $X_a$,
but there is no reason to beleive that sets $Y_{a,\gamma}$ are
closed. Note that there was a similar obstacle in
[Sh 351]. Similarly to [Sh 351], we overcome this by defining inductively the following sets $X_a^\oplus$ and $Z_{a,\gamma}$.

For simplicity in notation, let us introduce

\medskip

{\bf Definition 1.1.a.} (1)
Recall that a set $A$ of ordinals is said to be $\omega$-{\it closed\/}
if
$$\delta\in {\rm cl}(A)\,\,\&\,\,{\rm cf}(\delta)=\aleph_0
\implies \delta\in A.$$
We use ${\rm cl}$ to denote the ordinal closure.

(2) For a $\gamma\in
[\mu,\lambda)$ with ${\rm cf}(\gamma)>\aleph_0$, a  club
$C$ of $\gamma$ is $a$-{\it good\/} if 
$$\beta_1<\beta_2\in C\,\,\&\,\,{\rm cf}(\beta_1)
={\rm cf}(\beta_2)=\aleph_0\implies 
Y_{a,\beta_2}\cap  \beta_1 =Y_{a,\beta_1},$$
and
$$\beta\in C\,\,\&\,\,{\rm cf}(\beta)=\aleph_0
\implies \beta\in X_a\hbox{ and } Y_{a,\gamma}\hbox{ is }
\omega\hbox{-closed.}$$

\medskip

{\bf Remark.} 
Of course, we could
without loss of generality assume that our language has a constant for $\mu$
and so keep $\mu$ out of the definition
of $X_a$ and $Y_{a,\gamma}$. We may
skip the $\mu$ from similar definitions later.

\medskip

We define inductively
$$\eqalign{
X_a^\oplus\deq\{\gamma\in[\mu,\lambda):\,
& {\rm cf}(\gamma)=\aleph_0\,\,\&\,\,\gamma\in X_a
\,\,\&\,\, Y_{a,\gamma}\hbox{ is }\omega\hbox{-closed}\cr
&\vee  \kappa\ge{\rm cf}(\gamma)>\aleph_0
\,\,\&\,\,\hbox{there is a }a\hbox{-good club }C
\subseteq \gamma\cr
&\,\,\&\,\,
({\rm cf}(\delta)=\aleph_0\,\,\&\,\,\delta\in C)
\implies\delta\in X_a^\oplus\cr
&\vee {\rm cf}(\gamma)=1 \,\,\&\,\,
\hbox{there is a limit }\delta>\gamma
\hbox{ with }\delta\in X_a^{\oplus}\}.\cr}$$

For $\gamma\in X_a^\oplus$ we define inductively
$$Z_{a,\gamma}\deq\cases{
{\rm cl}(Y_{a,\gamma}\setminus\mu)\cap\gamma
&if ${\rm cf}(\gamma)=\aleph_0$\cr
\bigcap_{C \,a{\rm -good}\,\, {\rm club}\,\, {\rm of}
\,\,\gamma}
\bigcup_{\beta\in C\,\,\&\,\,{\rm cf}(\beta)=\aleph_0} Z_{a,\beta}
&if $\kappa\ge {\rm cf}(\gamma)>\aleph_0$ and $\gamma\in X_a^\oplus$\cr
Z_{a,\delta}\cap \gamma
&if ${\rm cf}(\gamma)=1$ and $\delta\in X_a^\oplus$\cr
& is
the minimal such limit $>\gamma$.
\cr}$$

We show that $\langle Z_{a,\gamma}:\,\gamma\in X_a^\oplus\rangle$ is
a square sequence on $X_a^\oplus$.
As $a$ is going to be fixed for some time, we may slip and say "good"
rather than $a$-good in the following.
Although we for most of the argument scholastically keep the $\bigcap$
over the good clubs of $\gamma$ in the definition of $Z_{a,\gamma}$
for $\gamma$ of uncountable cofinality, we invite the reader to check
that any two good clubs of $\gamma$ give the same value
to $\bigcup$ of the relevant $Z_{a,\beta}$. Hence we are in no
danger of intersecting too many sets. This argument in particular
shows that $Z_{a,\gamma}$ for such $\gamma$ is closed
and unbounded in $\gamma$.
Also note that in our definition of squares, successor ordinals
play no role, so the decision of what to put as a $Z_{a,\gamma}$ for
a successor $\gamma$ is quite arbitrary.

\medskip

{\bf Fact 1.1.b.} Suppose that $\gamma\in X_a^\oplus$ and $\beta\in Z_{a,\gamma}$. {\it Then\/}:

\item{(1)} $\beta\in X_a^\oplus$.
\item{(2)} If $\beta\in {\rm acc}(Z_{a,\gamma})$,
then $Z_{a,\beta}= Z_{a,\gamma}\cap\beta$.
\item{(3)} $\gamma$ a limit ordinal $\implies\sup (Z_{a,\gamma})=\gamma.$
\item{(4)} ${\rm otp}(Z_{a,\gamma})<\kappa^+$.
\item{(5)} $Z_{a,\gamma}$ is closed.

{\bf Proof.} (1)+(2) We prove the fist two items together
by induction on $\gamma$, dividing the discussion
into several cases.

\smallskip

$\underline{\hbox{Case I.}}$ ${\rm cf}(\beta)={\rm cf}(\gamma)
=\aleph_0$.

(1)--(2) Since $\beta\in Z_{a,\gamma}={\rm cl}(Y_{a,\gamma}\setminus\mu)\cap\gamma$, and
${\rm cf}(\beta)=\aleph_0$, by the $\omega$-closure of
$Y_{a,\gamma}$, we conclude that $\beta\in Y_{a,\gamma}$.
So, ${\rm Sk}_{\frak A}(a\cup\{\mu,\beta\})\subseteq 
{\rm Sk}_{\frak A}(a\cup\{\mu,
 \gamma\})$, hence
$Y_{a,\beta}\subseteq Y_{a,\gamma}\cap\beta$, and also
$$a\subseteq {\rm Sk}_{\frak A}(a\cup\{\mu,\beta\})\cap\mu
\subseteq {\rm Sk}_{\frak A}(a\cup\{\mu,\gamma\})\cap\mu
=a.$$
So, $\beta\in X_a$.
We now show $Y_{a,\beta}=Y_{a,\gamma}\cap\beta$,
from which it also follows that $Y_{a,\beta}$ is $\omega$-closed, hence $\beta\in X_a^\oplus$.

We already know that  $Y_{a,\beta}
\subseteq Y_{a,\gamma}\cap\beta$.
Now we proceed as in [Sh 430 1.1].  In $\frak A$
we can define just from $\mu,\beta$ a 1-1 onto function $f:\,\mu\into\beta$, as $\beta\in [\mu,\mu^+)$. The
$\elementary^{\ast}$-first such function, 
say $f^\ast$ is in ${\rm Sk}_{\frak A}(a\cup\{\mu,\beta\})$,
and also in ${\rm Sk}_{\frak A}(a\cup\{\mu,\gamma\})$, since this set contains $\beta$. In ${\rm Sk}_{\frak A}(a\cup\{\mu,\gamma\})$ ,
this function is a 1-1 onto function from $a$ onto
${\rm Sk}_{\frak A}(a\cup\{\mu,\gamma\})\cap\beta$.
In ${\rm Sk}_{\frak A}(a\cup\{\mu,\beta\})$, the range of this function is $Y_{a,\beta}$. Since $f^{\ast}$ is a fixed function, we
conclude that $Y_{a,\beta}=Y_{a,\gamma}\cap\beta$.

So,
$$Z_{a,\beta}=
{\rm cl}(Y_{a,\beta}\setminus\mu)\cap \beta=
{\rm cl}(Y_{a,\gamma}\cap\beta\setminus\mu)\cap\beta
={\rm cl}(Y_{a,\gamma}\setminus\mu)\cap\beta=
Z_{a,\gamma}\cap\beta.$$

$\underline{\hbox{Case II.}}$ ${\rm cf}(\beta)=\aleph_0\,\,\&\,\,{\rm cf}(\gamma)\in (\aleph_0,\kappa]$.

(1) Let $C\subseteq \gamma$ be an $a$-good club of $\gamma$.
By the definition of $Z_{a,\gamma}$, there is a $\delta\in C$ with
${\rm cf}(\delta)=\aleph_0$ such that $\beta\in Z_{a,\delta}$.
By the first case, $\beta\in X^\oplus_a$.

(2) For any $a$-good club $C$ of $\gamma$, let us denote by
$\delta_C$ the minimal element $\delta$ of $C$ with ${\rm cf}(\delta)=
\aleph_0$ such that $\beta\in Z_{a,\delta}$. Note that $\delta_C$
is well defined (by the definition of $Z_{a,\gamma}$), and $\delta_C
>\beta$. Then
$$Z_{a,\gamma}\cap\beta=\bigcap_{C \,a{\rm -good}\,\,
{\rm club}\subseteq
\gamma} \bigl [\bigcup_{\delta\in C\cap\delta_C\,\,\&\,\,{\rm cf}(\delta)=\aleph_0} (Z_{a,\delta}\cap\beta)
\cup\bigcup_{\delta\in C\setminus\delta_C\,\,\&\,\,{\rm cf}(\delta)=\aleph_0} (Z_{a,\delta}\cap\beta)\bigr]=$$
$$=\bigcap_{C \,a{\rm -good}\,\,{\rm club}\subseteq
\gamma}\bigl [\bigcup_{\delta\in C\cap\delta_C\,\,\&\,\,{\rm cf}(\delta)=\aleph_0} (Z_{a,\delta_C}\cap\delta\cap\beta)
\cup\bigcup_{\delta\in C\setminus\delta_C\,\,\&\,\,{\rm cf}(\delta)=\aleph_0} (Z_{a,\delta}\cap\delta_C\cap\beta)\bigr]=$$
$$=\bigcap_{C \,a {\rm-good}\,\,{\rm club}\subseteq
\gamma}\bigl[\bigcup_{\delta\in C\cap\delta_C\,\,\&\,\,{\rm cf}(\delta)=\aleph_0} Z_{a,\beta}\cap\delta\bigcup Z_{a,\beta}\bigr]=Z_{a,\beta}.$$

$\underline{\hbox {Case III.}}$ ${\rm cf}(\beta)\in(\aleph_0,\kappa]$ and
${\rm cf}(\gamma)=\aleph_0$.

(1) Then $\beta\in {\rm cl}(Y_{a,\gamma}\setminus\mu)\cap\gamma$. Let $A$ be an unbounded subset of $\beta$ with $A\subseteq Y_{a,\gamma}\setminus\mu$, and let $C$ be the ordinal closure
of $A$ in $\beta$. Hence $C$ is a club of $\beta$ and
$$\delta\in C\,\,\&\,\,{\rm cf}(\delta)=\aleph_0 \implies\delta\in
Y_{a,\gamma}\setminus\mu.$$
Suppose $\beta_1<\beta_2\in C$ and ${\rm cf}(\beta_1)={\rm cf}(\beta_2)=\aleph_0$. By the first case, we know that $\beta_1, 
\beta_2\in X_a^\oplus$ and
$$Y_{a,\beta_1}=Y_{a,\gamma}\cap\beta_1=Y_{a,\gamma}\cap\beta_2\cap\beta_1=Y_{a,\beta_2}\cap\beta_1.$$
Hence $C$ is a good club of $\beta$ and $\beta\in X_a^\oplus$.

(2) Suppose that $E$ is a good club of $\beta$ and $\delta\in E$
has cofinality $\aleph_0$. Without loss of generality, $E\subseteq C$,
where $C$ is as in the proof of (1) above. Hence $\delta\in Z_{a,\gamma}$, so by the first case, $Z_{a,\delta}=Z_{a,\gamma}\cap\delta$. Hence
$\bigcup_{\delta\in E\,\,\&\,\,{\rm cf}(\delta)=\aleph_0}Z_{a,\delta}=Z_{a,\gamma}\cap\beta$.
Hence $Z_{a,\beta}=Z_{a,\gamma}\cap\beta$.

$\underline{\hbox{Case IV.}}$
 ${\rm cf}(\beta),{\rm cf}(\gamma)\in (\aleph_0,\kappa]$.

(1) Let $C$ be a good club of $\gamma$, then there is a $\delta\in
C$ with ${\rm cf}(\delta)=\aleph_0$ and $\beta\in Z_{a,\delta}$.
By case III, $\beta\in X_a^\oplus$.

(2) Suppose $E$ is a good club of
$\gamma$. Let $\delta_E\in E$ be such that
$\beta\in Z_{a,\delta_E}$ and ${\rm  cf}(\delta_E)=\aleph_0$. Then,
as in Case III, we can find an $a$-good club $C$ of $\beta$ such that
$(\epsilon\in C\,\,\&\,\,{\rm cf}(\epsilon)=\aleph_0)\implies
\epsilon\in Y_{a,\gamma}$.
In particular (by the first case) for any $\epsilon\in C$
with ${\rm cf}(\epsilon)=\aleph_0$, we have
$Z_{a,\epsilon}=Z_{a,\delta_E}\cap\epsilon$.
Hence $\bigcup \{ Z_{a,\epsilon}
:\,\epsilon\in C\,\,\&\,\,{\rm cf}(\epsilon)=\aleph_0\}
=Z_{a,\delta_E}\cap\beta$. By the third case,
this is equal to $Z_{a,\gamma}\cap\beta$.
In our calculation of $Z_{a,\beta}$ we can without loss of
generality restrict ourselves to the good clubs of
$\beta$ which are subsets of $C$.
Hence $Z_{a,\beta}=Z_{a,\gamma}\cap\beta$.

$\underline{\hbox{Case V.}}$ ${\rm cf}(\beta)=1$ and ${\rm cf}(\gamma)\in [\aleph_0,\kappa]$.

(1) By the definition of $X_a^\oplus$, we have that $\beta\in X_a^\oplus$.

(2) Does not apply.

$\underline{\hbox{Case VI.}}$ ${\rm cf}(\beta)={\rm cf}(\gamma)=1.$

(1) Let $\delta$ be a limit ordinal $>\gamma$ such that
$\delta\in X_a^\oplus$. Then $\beta<\delta$, so $\beta\in X_a^\oplus$.

(2) Does not apply.

$\underline{\hbox{Case VII.}}$ ${\rm cf}(\beta)>1$ and ${\rm cf}(\gamma)=1.$

(1) Let $\delta>\gamma$ be a limit ordinal such that $Z_{a,\gamma}=Z_{a,\delta}\cap\gamma$, and use case V.

(2) If $\beta\in{\rm acc}(Z_\gamma)$, then
$\beta\in{\rm acc}(Z_\delta)$, so $Z_\beta=Z_\delta\cap\beta=
Z_\gamma\cap\beta$, by previous cases.

\smallskip

We proceed to the proof of (3)--(5).

(3) Suppose that $\gamma$ is a limit ordinal.
First suppose that ${\rm cf}(\gamma)=\aleph_0$.
In $\frak A$, there is a cofinal function $f:\,\aleph_0
\into\gamma$, definable from $\gamma$ only. Hence the
first such $f$ is an element of $Y_{a,\gamma}$. Then
${\rm Ran}(f)$ is an unbounded subset of $\gamma\cap
Y_{a,\gamma}$, so $\sup(Z_{a,\gamma})\ge\sup(Y_{a,\gamma}
\cap\gamma)=\gamma$.

If ${\rm cf}(\gamma)>\aleph_0$,
then we know that for every $\delta\in Z_{a,\gamma}$ with
${\rm cf}(\delta)=\aleph_0$ we have $\delta=\sup (Z_{a,\delta})$,
so the conclusion follows from the definition of $Z_{a,\gamma}$.

(4) This follows since $\card{Y_{a,\beta}}\le\kappa$ for any $\beta\in X_a$, so putting all the definitions together,
$\card{Z_{a,\gamma}}\le\kappa$.

(5) If ${\rm cf}(\gamma)=
\aleph+)$, this is implicitly stated in the definition.
Next, it is easy to check that for ${\rm cf}(\gamma)>
\aleph_0$
(see the paragraph before Fact 1.1.b).
Finally, for $\gamma$ a successor, we also
obtain a closed set via our definitionof $Z_{a,\gamma}$.\eop${}^{1.1.b.}$

\medskip

So we have now to prove that we can choose $\mu$ many $a$ such that all $\gamma\in [\mu,\lambda)$ of cofinality $\le\kappa$ are
in some $X_a^\oplus$.
We shall use the following Theorem of Saharon Shelah, which is a consequence
 of Theorem 1.4 of [Sh 580].  We note that in fact a stronger version of Theorem C follows from [Sh 580, 1.4], where $\mu$ is not required to be a singular $>\kappa$, but just to be above some
finitely many cardinal successors of $\kappa$.

\medskip

{\bf Theorem C (Shelah)} Suppose that $\kappa$
and $\mu$ are as above.
{\it Then\/} there is a $\PP\subseteq [\mu]^{\le\kappa}$ with
$\card{\PP}=\mu$, such that for every $\theta\ge\kappa$ and
$x\in H(\theta)$, we can find a continuously increasing sequence
$\bar N=\langle N_i:\,i<\kappa^+\rangle$ such that:
\item{-} $N_i\in N_{i+1}$ for $i\le\kappa^+$.
\item{-} $x\in N_i\elementary \langle H(\theta),\in,\elementary^\ast\rangle$ for all $i$.
\item{-} ${\card{N_i}}=\kappa$ for all $i$ and $\kappa\subseteq N_i$.
\item{-} For every club $E$ of $\kappa^+$, there is an $i\in E$ such that $N_i\cap\mu\in \PP$.

\medskip

Now we fix a $\PP\subseteq [\mu]^{\le\kappa}$ as in Theorem C,
so $\card{\PP}=\mu$.
We claim that $\bigcup\{X_a^\oplus:\,a\in \PP\}$ contains all
$\gamma\in[\mu,\lambda)$ with ${\rm cf}(\gamma)\le\kappa$.
This suffices, as we have just proved that
$\langle Z_{a,\gamma}:\,\gamma\in X^\oplus_a\rangle$ is a square
sequence, for any $a\in [\mu]^{\le\kappa}$. It suffices to prove the following 

\medskip

{\bf Claim 1.1.c.} For every $\gamma\in[\mu,\lambda)$ with $\kappa\ge{\rm cf}(\gamma)$, there is an $a\in\PP$ such that
$\gamma\in X_a^\oplus$.

{\bf Proof.} It suffices to prove this for $\gamma\in[\mu,\lambda)$ with $\kappa\ge{\rm cf}(\gamma)
\ge\aleph_0$.
Let us define for such $\gamma$,
$$\bar\gamma\deq\cases{\emptyset      &if cf$(\gamma)=\aleph_0$,\cr
\{\gamma_\epsilon:\,\epsilon<{\rm cf}(\gamma) \}   &if ${\rm cf}(\gamma)>\aleph_0$,\cr}
$$
 where $\{\gamma_\epsilon:\,\epsilon<{\rm cf}(\gamma)\}$ is an
increasing enumeration of a club of $\gamma$. Let $\theta$ be large
enough, so that $\frak A\in H(\theta)$, and let $x=\langle
\gamma,\bar{\gamma},\frak A\rangle$.
Let ${\frak B}\deq\langle H(\theta),\in,\elementary^\ast\rangle$.

Using Theorem C, we get a sequence $\langle N_i:\,i<\kappa^+\rangle$ and club $E$ of $\kappa^+$ which
exemplify the theorem for our chosen $x$.

Now, for any $i$, we know that $N_i\cap\gamma\in N_{i+1}$. 
Hence, ${\rm cl}(N_i\cap\gamma)\cap\gamma\subseteq N_{i+1}$.

So, if ${\rm cf}(i)>\aleph_0$, then $N_i\cap\gamma$ is $\omega$-closed. Note that for any $i$, we have
$\{\frak A\}\cup(N_i\cap \mu)\cup\{\gamma\}\subseteq N_i\elementary
{\frak B}$, so ${\rm Sk}_{\frak A}\bigl((N_i\cap \mu)\cup\{\gamma\}\bigr)\subseteq N_i$,
hence
${\rm Sk}_{\frak A}\bigl((N_i\cap \mu)\cup\{\gamma\}\bigr)
\cap \mu= N_i\cap\mu.$

So, if ${\rm cf}(\gamma)=\aleph_0$ and $i\in E$ is such that
${\rm cf}(i)>\aleph_0$, then $\gamma\in X_{N_i\cap\mu}$.
Similarly, if $a=N_i\cap\mu$, then $Y_{a,\gamma}=N_i\cap\gamma$.
\eop${}_{1.1.c}$

The conclusion of the theorem follows easily, the sets
with squares are $X_a^\oplus$ for $a\in\PP$.\eop${}_{1.1.}$

\medskip

{\bf Remark 1.2}(0) Note that we have obtained an alternative proof that
under the assumptions of 1.1, we have $S^\lambda_{\le\kappa}\in I[\lambda]$.
This was proved in [Sh 420, 2.8].

(1) Theorem 1.1. strengthens [Sh 237e, 2]
for $\mu$ singular, as [Sh 237e,2] had the same conclusion under
$\mu^{\le\kappa}=\mu$
instead of $\cf\bigl([\mu]^{\le\kappa},\subseteq\bigr)=\mu$.

(2) If $\aleph_0<\cf (\kappa)=\kappa\le\cf (\mu)<\mu$, what is the strength
of the assumption $\cf ([\mu]^{\le\kappa},\subseteq)=\mu$? In [Sh 430, 1.3.]
it is proved that this follows from ${\rm pp}(\mu)=\mu^+$.
If $\cf ([\mu]^{\aleph_0},\subseteq)=\mu$ and for all $\theta\in(\kappa,\mu)$
we have $\theta>\cf(\theta)<\cf (\mu)\implies {\rm pp}(\theta)\le\mu$, then
$\cf ([\mu]^{\le\kappa},\subseteq)=\mu$.

A particular situation in which our theorem applies and is not implied by the
previously known theorems, is for example $\mu=\aleph_{\omega_{17}}$ and
$\kappa=\aleph_{13}$ (see [Sh 400, Why the HELL is it four?]).

\medskip

{\it 1.3. Acknowledgement\/}
We would like to thank  James Cummings for noticing a mistake in an
earlier version of the theorem.

\medskip

We shall now turn our attention to successors of singular strong limits,
for which we can prove that a weak version of the square property
holds. It will be useful to define the notion of a {\it square pretender\/},
as follows.

\medskip

{\bf Definition 1.4.} Suppose that $\kappa<\mu$ are cardinals, $\kappa$
is regular,
and $e\subseteq\mu$. A {\it square pretender on\/} $e$
{\it of depth\/} $\kappa$ is
a sequence
$$\bigl\langle S_i,\bar{d}^i=\langle d_\gamma^i:\,\gamma\in S_i\rangle,
\bar{s}^{\gamma,i}=\langle
\alpha_\zeta^{\gamma,i}:\,\zeta\in d_\gamma^i\rangle:\,i<\kappa\bigr\rangle$$
such that:

\item{$(a)$} $\cup_{i<\kappa} S_i\supseteq\{\beta\in e:\,\cf(\beta)\neq
\kappa\}$ and
$\cup_{i<\kappa}\cup_{\gamma\in S_i}
\{\alpha_\zeta^{\gamma,i}:\,\zeta\in d_\gamma^i\}\supseteq
\{\beta\in e:\,\cf(\beta)\neq\kappa\}.$
\item{$(b)$} $\langle \alpha_\zeta^{\gamma,i}:\,\zeta\in d_\gamma^i\rangle$
is an increasing sequence of elements of $\gamma$.
\item{$(c)$} If $\zeta\in d_\gamma^i$, then $\alpha_\zeta^{\gamma,i}\in S_i$
and $d_{\alpha_\zeta^{\gamma,i}}^i\subseteq d_\gamma^i$ and
$\bar{s}^{\alpha_\zeta^{\gamma,i},i}=\langle \alpha_{\zeta^{'}}^{\gamma,i}:\,
{\zeta^{'}}\in d_{\alpha^{\gamma,i}_\zeta}^i\rangle.$

\medskip

Before we state the following theorem, we remind the reader of the
following.

We shall be concerned with $\lambda=\mu^+$ for $\mu$ a strong limit
singular cardinal. By [Sh 108] or [Sh 88a] in this situation there is
a maximal $W^\ast$ in $I[\lambda]$,
which is unique up to a nonstationary set. In other words, for every
$W\subseteq\lambda$, we have that $W\in I[\lambda]\iff
W\setminus W^\ast\in NS[\lambda]$. (A reader unfamiliar with the
ideal $I[\lambda]$ will find a lot about it in \S3,
including the meaning of $\langle C_\delta:\,\delta<\lambda\rangle$
witnessing that $W^\ast\in I[\lambda]$, which is needed for the
statement of the following theorem.)

On to the theorem:

\medskip

{\bf Theorem 1.5.} Suppose that $\mu$ is a strong limit singular cardinal
of cofinality $\kappa$ and $\lambda=\mu^+$. Let $W^\ast$ be the maximal
element of $I[\lambda]$ and let $\langle C_\delta:\,
\delta<\lambda\rangle$ witness that $W^\ast\in I[\lambda]$. Let $E
=\acc (E_0)$, where $E_0$ is a club
of $\lambda$ such that for every $\delta\in W^\ast\cap E_0$, we have
that $C_\delta$ is a club of $\delta$ of order type $\cf(\delta)<\delta$
(see 3.0.2).

{\it Then,\/} there is a sequence
$$\bigl\langle S_j, \bar{d}^j=\langle d_\gamma^j:\,\gamma\in S_j\rangle,
\bar{s}^{\gamma,j}=\langle \alpha_\zeta^{\gamma,j}:\,\zeta\in d_\gamma^j
\rangle:\,j<\mu\bigr\rangle$$
such that for every $\delta\in W^\ast\cap E$ and for every club $e
\subseteq\acc(C_\delta)$, there is a sequence $\langle j_i:\,i<\kappa
\rangle$ in $\mu$ such that
$$\langle S_{j_i},\bar{d}^{j_i},\bar{s}^{\gamma,j_i}:\,i<\kappa\rangle$$
is a square pretender on $e$ of depth $\kappa$.

{\bf Proof.} Let us fix an increasing sequence $\langle\mu_i:\,
i<\kappa\rangle$ of cardinals such that $\mu=\Sigma_{i<\kappa}
\mu_i$. We also choose
by induction on $\alpha<\lambda$, sets $a^\alpha_i$ for
$i<\kappa$ with the following properties:

\item{$(a)$} $\alpha=\cup_{i<\kappa}a^\alpha_i$.
\item{$(b)$} $\card{a^\alpha_i}\le\mu_i.$
\item{$(c)$} $\alpha\in a^\beta_i\implies a^\alpha_i\subseteq
a^\beta_i$.
\item{$(d)$} $i<j<\kappa\implies a^\alpha_i\subseteq
a^\alpha_j.$
\item{$(e)$} $\card{C_\alpha}\le\mu_i\implies C_\alpha\subseteq a^\alpha_i$.

Let us also define for $\alpha<\beta<\lambda$,
$$c(\alpha,\beta)\deq\min\{i<\kappa:\,\alpha\in a^\beta_i\},$$
so
$\alpha<\beta<\gamma\implies c(\alpha,\gamma)\le\max\bigl\{
c(\alpha,\beta),c(\beta,\gamma)\bigr\}$, by $(c)$ above.

Now let us fix a $\delta\in W^\ast\cap E$ and let $\theta=\cf(\delta)$,
so $\theta$ is a regular cardinal $<\mu$. Suppose that
$e\subseteq\acc(C_\delta)$ is a fixed club of $\delta$. Therefore $\otp(e)
=\theta$ (since $\otp (C_\delta)=\cf(\delta))$. We define for all $i<\kappa$
such that $\mu_i\ge\theta$,
$$A_{\delta,i}=A_{\delta,i,e}\deq
\bigl\{\alpha\in e:\,a^\alpha_i\cap\nacc(C_\delta)\hbox{ is unbounded in }
\alpha\bigr\}.$$
If $\mu_i<\theta$, we define $A_{\delta,i}=A_{\delta,i,e}\deq\emptyset$.

We prove some facts about the just defined sets, which will prepare the
ground for further definitions.

\smallskip

{\bf Fact 1.5.a.} If $\alpha_1<\alpha_2\in A_{\delta,i}$, then
$a^{\alpha_1}_i$ is a bounded subset of $a^{\alpha_2}_i$
(hence $\otp(a^{\alpha_1}_i)$\break$<\otp(a^{\alpha_2}_i)$.)

{\bf Proof.} Since $\alpha_2\in A_{\delta,i}$, we can find
a $\beta\in\nacc(C_\delta)\cap a^{\alpha_2}_i$ which is $>\alpha_1$.
By $\beta\in\nacc(C_\delta)$, we have $C_\beta=C_\delta\cap\beta$
(see 3.0.2). Now, since $\alpha_1\in e$, in particular
$\alpha_1\in C_\delta$,
so $\alpha_1\in C_\beta$. By $(e)$ and $(c)$ we have
$a^{\alpha_1}_i\subseteq a^\beta_i$. But $a^\beta_i\subseteq a^{\alpha_2}_i$
as $\beta\in a^{\alpha_2}_i$. Obviously, by $\alpha_2\in A_{\delta,i}$, we have
$\sup(a^{\alpha_2}_i)=\alpha_2>\alpha_1\ge\sup(a^{\alpha_1}_i).$\eop${}_{1.5.a.}$

\smallskip

{\bf Fact 1.5.b.} $\langle A_{\delta,i}:\,i<\kappa\rangle$ is an increasing
sequence of subsets of $e$.

{\bf Proof.} This follows, since $a^\alpha_i$ are increasing.\eop${}_{1.5.b.}$

\smallskip

{\bf Fact 1.5.c.} $\cup_{i<\kappa} A_{\delta,i}\supseteq\{\gamma\in e:\,
\cf(\gamma)\neq\kappa\}$.

{\bf Proof.} For any $\gamma\in e\subseteq\acc(C_\delta)$, we have that
$\gamma=\sup(C_\delta\cap\gamma)=\sup\bigl(\nacc(C_\delta)\cap\gamma\bigr)$,
so $\gamma=\sup\bigl(\nacc(C_\delta)\cap\cup_{i<\kappa}a^\gamma_i\bigr).$
If $\cf(\gamma)\neq\kappa$, then there is an $i<\kappa$ such that
$\gamma=\sup\bigl(\nacc(C_\delta)\cap a^\gamma_i\bigr)$.
As $a^\gamma_i$ are increasing with $i$,
there is such $i$ such that $\mu_i\ge \theta$.\eop${}_{1.5.c.}$

\smallskip

{\bf Remark 1.5.d.}
{\it (this remark is not used later in the proof)\/} Suppose $e\subseteq E_0$.
If $\gamma\in e$ is such that $A_{\delta,i}\cap\gamma$ is
stationary in $\gamma$
and $\card{C_\gamma}\le\mu_i$, then $\gamma\in A_{\delta,i}$.

{\bf Proof.} Since $\gamma\in e\subseteq \acc(C_\delta)$, the set
$T\deq A_{\delta,i}\cap C_\delta\cap C_\gamma$
is stationary in $\gamma$
(as $\gamma\in E_0$). If $\beta\in T$, by $(e)$ in the definition
of $a's$, we have $a^\beta_i\subseteq a^\gamma_i$. Therefore
$$\sup\bigl(a^\gamma_i\cap\nacc(C_\delta)\bigr)\ge
\sup_{\beta\in T} \bigl(a^\beta_i\cap\nacc(C_\delta)\bigr)=
\sup_{\beta\in T}\beta=\gamma.
\hbox{\eop}{}_{1.5.d.}$$

\smallskip

{\it Continuation of the proof of 1.5.\/} Let us now fix an
$i(\ast)<\kappa$.

We enumerate $A_{\delta,i(\ast)}$ increasingly as
$$A_{\delta,i(\ast)}=\{\alpha^\delta_\epsilon:\,\epsilon <
\otp(A_{\delta,i(\ast)})\le\theta\},$$
and set $\alpha^\delta_\theta=\delta$. For $\epsilon\le\theta$
for which $\alpha^\delta_\epsilon$ is defined, we define
$$d^{i(\ast)}_\epsilon\deq\{\zeta<\epsilon:\,
c(\alpha^\delta_\zeta,\alpha^\delta_\epsilon)\le i(\ast)\}$$
and for $\zeta\le\epsilon\le\theta$
$$b^{i(\ast)}_{\zeta,\epsilon}\deq\otp(\alpha^\delta_\zeta\cap
a^{\alpha^\delta_\epsilon}_{i(\ast)}).$$
We define a partial function $\epsilon:\lambda\to\theta$ by setting
$\epsilon(\alpha)$ to be the (unique by 1.5.f) $\epsilon$ such that we can
find a sequence $\langle\alpha_\zeta:\,\zeta\in d^{i(\ast)}_\epsilon
\rangle$ in $\alpha$ with the following properties:

\item{$(A)$} $\langle \alpha_\zeta:\,\zeta\in d^{i(\ast)}_{\epsilon}
\rangle$ increases with $\zeta$.
\item{$(B)$} $c(\alpha_\zeta,\alpha)\le i(\ast).$
\item{$(C)$} $\otp(\alpha_\zeta\cap a^\alpha_{i(\ast)})=b^{i(\ast)}_{\zeta,
\epsilon}.$
\item{$(D)$} If $\zeta_1<\zeta_2$ are in $d^{i(\ast)}_\epsilon$ and
$c(\alpha^\delta_{\zeta_1},\alpha^\delta_{\zeta_2})\le i(\ast)$,
{\it then,\/}
$$c(\alpha_{\zeta_1},\alpha_{\zeta_2})=c(\alpha^\delta_{\zeta_1},
\alpha^\delta_{\zeta_2})$$
and
$$\otp(\alpha_{\zeta_1}\cap a^{\alpha_{\zeta_2}}_{i(\ast)})
=\otp(\alpha^\delta_{\zeta_1}\cap a^{\alpha^\delta_{\zeta_2}}_{i(\ast)}).$$
\item{$(E)$} $\otp(a^\alpha_{i(\ast)})=\otp(a^{\alpha^\delta_\epsilon}_{
i(\ast)})$ by an order preserving isomorphism which exemplifies that
$\otp(a^{\alpha_\zeta}_{i(\ast)})=\otp(a^{\alpha^\delta_{\zeta}}_{i(\ast)})$
for all $\zeta\in d^{i(\ast)}_{\epsilon}.$

We prove several facts about the partial function $\epsilon$
and sets $d_\epsilon^{i(\ast)}$.

\smallskip

{\bf Fact 1.5.e.} If $\epsilon_1\in d^{i(\ast)}_\epsilon$, then
$d^{i(\ast)}_{\epsilon_1}\subseteq d^{i(\ast)}_\epsilon$.

{\bf Proof.} If $\zeta\in d^{i(\ast)}_{\epsilon_1}$, then $\zeta<\epsilon_1$
and $c(\alpha^\delta_\zeta,\alpha^\delta_{\epsilon_1})\le
i(\ast)$. Since $\epsilon_1\in d^{i(\ast)}_\epsilon$, also $\epsilon_1<\epsilon$
and
$c(\alpha^\delta_{\epsilon_1}, \alpha^\delta_\epsilon)\le i(\ast)$. So
$\zeta<\epsilon$ and $c(\alpha^\delta_\zeta,\alpha^\delta_\epsilon)
\le i(\ast).$\eop${}_{1.5.e.}$

\smallskip

{\bf Fact 1.5.f.} If $\alpha<\lambda$, there is at most one
$\epsilon$ and $\langle\alpha_\zeta:\,\zeta\in d^{i(\ast)}_\epsilon\rangle$
which satisfy $(A)$--$(E)$.

{\bf Proof.} Suppose first that $\langle \alpha_\zeta:\,\zeta\in 
d^{i(\ast)}_\epsilon\rangle$ and $\langle \alpha^{'}_\zeta:\,
\zeta\in d^{i(\ast)}_\epsilon\rangle$ both exemplify
$(A)$--$(E)$. Then for each $\zeta\in d^{i(\ast)}_\epsilon$, by $(B)$
above
$$\alpha_\zeta,\alpha^{'}_\zeta\in a^\alpha_{i(\ast)}.$$
By $(C)$ above
$$\otp(\alpha_\zeta\cap a^\alpha_{i(\ast)})=\otp(\alpha^{'}_\zeta\cap
a^\alpha_{i(\ast)}),$$
so it must be that $\alpha_\zeta=\alpha_{\zeta}^{'}$, for all $\zeta\in
d^{i(\ast)}_\epsilon.$

Suppose then that $\langle \alpha_\zeta:\,\zeta\in d^{i(\ast)}_{\epsilon_1}
\rangle$ and $\langle\beta_\zeta:\,\zeta\in d^{i(\ast)}_{\epsilon_2}
\rangle$ both exemplify $(A)$--$(E)$ for $\epsilon_1<\epsilon_2$.
By $(E)$
$$\otp(a^{\alpha^\delta_{\epsilon_1}}_{i(\ast)})=
\otp(a^\alpha_{i(\ast)})=\otp(a^{\alpha^\delta_{\epsilon_2}}_{i(\ast)}),$$
which is a contradiction with $\alpha^\delta_{\epsilon_1}<\alpha^\delta_{
\epsilon_2}\in A_{\delta,i(\ast)}$, by 1.5.a.\eop${}_{1.5.f.}$

\smallskip

{\bf Fact 1.5.g.} If $\epsilon(\alpha)$ is well defined, and witnessed
by $\langle\alpha_\zeta:\,\zeta\in d^{i(\ast)}_{\epsilon(\alpha)}
\rangle$, {\it then,\/} for every $\zeta\in d^{i(\ast)}_{\epsilon(\alpha)}$,
we have that $\epsilon(\alpha_\zeta)=\zeta$ and this is witnessed by
$\langle \alpha_{\zeta^{'}}:\,\zeta^{'}\in d^{i(\ast)}_\zeta\rangle$.

{\bf Proof.} By Fact 1.5.e, we have that $d^{i(\ast)}_\zeta\subseteq
d^{i(\ast)}_{\epsilon(\alpha)}$, so $\langle\alpha_{\zeta^{'}}:\,
\zeta^{'}\in d^{i(\ast)}_{\zeta}\rangle$ is well defined. We need to check that
$(A)$--$(E)$ are satisfied.

$(A)$ is obviously true, so consider $(B)$. By $(D)$ for $\epsilon(\alpha)$
and the definition of $d^{i(\ast)}_\zeta$,
if $\zeta^{'}\in d^{i(\ast)}_\zeta$,
then $c(\alpha_{\zeta^{'}},\alpha_\zeta)=
c(\alpha^\delta_{\zeta^{'}},\alpha^\delta_\zeta)\le i(\ast).$

To see $(C)$, use again $(D)$, so $\otp(\alpha_{\zeta^{'}}\cap a^{\alpha_\zeta}_
{i(\ast)})=\otp(\alpha^\delta_{\zeta^{'}}\cap a^{\alpha^\delta_\zeta}_{i(\ast)})$,
which is by definition equal to $b^{i(\ast)}_{\zeta^{'},\zeta}$.

Now, $(D)$ follows from $(D)$ for $\epsilon(\alpha)$, and the fact that
$d_\zeta^{i(\ast)}\subseteq d^{i(\ast)}_{\epsilon(\alpha)}$. The last
statement also implies $(E)$.\eop${}_{1.5.g.}$

\smallskip

{\it Continuation of the proof of 1.5.} Now we can set
$$S_{i(\ast)}=S_{\delta,i(\ast),e}\deq
\{\alpha\in\lambda:\,\epsilon(\alpha)\hbox{ is well defined }\}.$$
Note that $S_{i(\ast)}\supseteq A_{\delta,i(\ast)}$, as for
$\alpha=\alpha^\delta_{\epsilon(\ast)}$ we have $\epsilon(\alpha)=
\epsilon(\ast).$
Also note that $\cup_{i<\kappa}\cup_{\alpha\in S_i}
\{\alpha_\zeta:\,\zeta\in d^i_{\epsilon(\alpha)}\}\supseteq
\cup_{i<\kappa} A_{\delta,i}$.

Therefore $\bigl\langle S_i,\langle d^i_{\epsilon(\alpha)}:\,\alpha
\in S_i\rangle, \langle \alpha_\zeta:\,\zeta\in d^i_{\epsilon(\alpha)}
\rangle:\,i<\kappa\bigr\rangle$ is a square pretender of depth
$\kappa$ on $e$.

Now notice that the choice of this square pretender only depended on the
following parameters:

\item{{\it (i)\/}} $\theta<\mu$.
\item{{\it (ii)\/}} $\langle \otp(A_{\delta,i}):\, i<\kappa\rangle$, which
is a $\kappa$-sequence of elements of $\theta+1$.
\item{{\it (iii)\/}} $\bigl\langle\langle d_\epsilon^i:\,\epsilon\le\theta
\rangle:\,i<\kappa\bigr\rangle$, and each $d^i_\epsilon$ is a subset of
$\epsilon$.
\item{{\it (iv)\/}} $\bigl\langle\langle b^i_{\zeta,\epsilon}:\,\zeta\le\epsilon
\le\theta\rangle:\,i<\kappa\rangle$ and each $b^i_{\zeta,\epsilon}$ is an
element of $\mu^+_i$.

By the fact that $\mu$ is a strong limit, all possible such choices can be
enumerated in type $\mu$, so the theorem is proved.\eop${}_{1.5.}$
\bigskip

{\bf \S2. On singular ex-inaccessible and outside guessing of clubs.}
In this section we consider the following question. Suppose we
start with two universes $V_1\subseteq V_2$
of set theory and a cardinal $\kappa$ which
is inaccessible in $V_1$. We assume $\kappa$ is singular
in $V_2$, and then
want to ``guess'' clubs of $\kappa$ from $V_1$
(we refer to this as to ``outside guessing of clubs").
The history of this question and related results of Shelah, and Gitik,
were explained in the introduction.
 
We have two ways of ``guessing'':
the first one is to find an unbounded subset of 
(or, equivalently, find a club) of $\kappa$ in $V_2$
which is almost contained in every club of $\kappa$ in $V_1$.
The other ``guessing'' is by a proper filter.
To obtain these guessings we need additional assumptions, which
go in two directions. One is the cardinal arithmetic in $V_1$.
If the extension $V_2$ was obtained through forcing, the other set of
restrictions can be regarded as
speaking on the chain condition of the forcing used.
In fact, these restrictions are about certain covering numbers.

Throughout this section, if
we are simultaneously speaking of two universes of set theory,
$V_1$ and $V_2$, such that $V_1\subseteq V_2$, and we have not
specified in which one we carry the argument, then we mean $V_2$.
The symbol $\Delta$ stands for the diagonal intersection.

We now proceed to the results.

\medskip

{\bf Theorem 2.0.} Assume that:

\item{{\it (i)}\/} $V_1\subseteq V_2$ are transitive classes
containing the
ordinals
and modeling $ZFC$.
\item{{\it (ii)\/}} $V_1\satisfies\,``\kappa\hbox{ is inaccessible.}"$
\item{{\it (iii)\/}}
$V_2\satisfies \,``\kappa\hbox{ is a singular cardinal},
\cf(\kappa)=\theta$ and $ \kappa^+=(\kappa^+)^{V_1}".$
\item{{\it (iv)\/}} $ V_1\satisfies\,``\cf\bigl(\Club
 (\kappa),\supseteq\bigr)                                   
 \leq\kappa^+".$

{\it Then\/}:

\item{(1)} In $V_2$, we can find an
unbounded $C^\ast\subseteq \kappa$
such that:
 $$E\in {\Club^{V_1}}(\kappa)\implies C^\ast\setminus E\hbox{
 is bounded.}$$
\item{(2)} If $\gamma<\kappa$ is fixed, we can also demand that
 $$\beta\in\nacc(C^\ast)\implies{\cf}^{V_2}(\beta)>\gamma.$$
 \item{(3)} If $V_2\satisfies\,``\theta>\aleph_0"$, we can add in (1):
 $$\alpha\in\nacc(C^\ast)\implies{\cf}^{V_2}(\alpha)\ge\sup(C^\ast\cap\alpha)
 .$$

{\bf Proof.} (1) In $V_1$,
by
[Sh 351, \S4]=[Sh 365, 2.14],
we can
 find sets $S_i\,(i<\kappa)$ such that:

 \item{$(\alpha)$} $\cup_{i<\kappa}
 S_i=\{\delta<\kappa^+:\,\cf(\delta)<\kappa\}.$

 \item{$(\beta)$}  For each $i<\kappa$, there is a square $\langle
 C_\alpha^i:
 \,\alpha\in S_i\rangle$ on $S_i$ (see 1.0 for the definition).

In addition to $(\alpha)$ and $(\beta)$,
 since $\kappa$ is a limit cardinal, (again by [Sh 351, \S 4]=[Sh 365, 2.14])
we can assume
that:

 \item{$(\gamma)$} For each $i<\kappa$ there is a $\mu_i<\kappa$,
 such that
 for all $\alpha\in S_i$
we have $\otp(C^i_\alpha)<\mu_i$.

Since $(\kappa^+)^{V_1}=(\kappa^+)^{V_2}$, and
it is easily seen that $(\beta)$ and $(\gamma)$ are
absolute
between $V_1$ and $V_2$,
the sets $S_i\,(i<\kappa)$ will still satisfy $(\beta)-(\gamma)$
in $V_2$, and, by $(\alpha)$,
we shall have $\cup_{i<\kappa}S_i\supseteq\{\delta
<\kappa^+:\,{\cf}^{V_2}(\delta)\neq\theta\}$.

Now we argue in $V_2$.
We fix a regular $\mu^\ast\in(\theta,\kappa)$.

Since $S_{\mu^{\ast}}^{\kappa^+}\deq\{\alpha<\kappa^+:\,\cf^{V_2}(\alpha)=
\mu^{\ast}\}$ is stationary in $\kappa^+$
and included in $\cup_{i<\kappa}S_i$, we can fix an $i(\ast)<\kappa$ such that
$S_{i(\ast)}\cap S_{\mu^{\ast}}^{\kappa^+}$ is stationary. If we take
any club $C$ of $\kappa^+$ which is in $V_2$, and if $\alpha\in
S_{i(\ast)}\cap S_{\mu^{\ast}}^{\kappa^+}\cap \acc(C)$, then
${\cf}^{V_2}(\alpha)=\mu^\ast$ and $C\cap C^{i(\ast)}_\alpha$ is a club of
$\alpha$.
To summarize, we can conclude

 \item{$(\delta)$} There is an
 $i(\ast)<\kappa$ such that, in $V_2$:

 {\parindent=35pt

 For every club $C\in V_2$ of $\kappa^+$, for
 stationarily
 many $\alpha\in S_{i(\ast)}$,
$$ C_\alpha^{i(\ast)}\cap C\hbox{ is a club of }\alpha\,\,\&\,\,
 {\cf}^{V_2}(\alpha)=\mu^\ast.$$
 }
{\it Observation.}

In $V_1$, by $\cf^{V_1}(\Club(\kappa),\supseteq)\leq\kappa^+$
we can fix a
 sequence $\langle E_\alpha:
 \,\alpha<\kappa^+\rangle\in V_1$ of clubs in $\kappa$, such that:

 \item{(i)} For every $E\in\Club^{V_1}(\kappa)$, for some $\alpha$
we have
 $E\supseteq
 E_\alpha$.
 \item{(ii)} $\beta<\alpha\implies E_\beta\supseteq^{\ast}E_\alpha$.
 \item{(iii)} $\beta\in C_\alpha^{i(\ast)}\implies E_\beta\supseteq
 E_\alpha.$

[Why?
This can be done by induction: suppose that $\langle D_\alpha:\,\alpha
<\kappa^+\rangle\in V_1$ is a sequence of clubs of $\kappa$ which is cofinal in
$({\Club}^{V_1}(\kappa),\supseteq)$. Let $E_0=D_0$, and
suppose that $\langle E_\beta:\,\beta
<\alpha\rangle$ are given for an $\alpha<\kappa^+$.
Note that the following is well defined
$$\gamma_\alpha\deq\min\{\gamma:\,(\forall\beta<\alpha)
\,(D_\gamma\subseteq^\ast E_\beta)\}.$$
Namely, $E^\beta
\deq\Delta_{\beta<\alpha} E_\beta\in \hbox{Club}^{V_1}(\kappa)$
and satisfies $\forall\beta<\alpha (E^\beta\subseteq^\ast E_\beta).$
So, there is some $\gamma$ such that $D_\gamma\subseteq^\ast E^\beta$.

Now, if $\gamma_\alpha\in S_{i(\ast)}$, let
$$E_\alpha\deq D_\alpha\cap D_{\gamma_\alpha}
\cap\bigcap_{\beta\in C^{i(\ast)}_\alpha}E_\beta.$$
If $\alpha\notin S_{i(\ast)}$, we simply skip the last term in the above
intersection. 
Note that we are using $(\gamma)$ above to assure that each $E_\alpha$
is a club.]

The rest of the argument takes place in $V_2$. We
let $\langle\xi_\epsilon:\,\epsilon<\theta\rangle$ be
 a strictly increasing cofinal sequence in $\kappa$,
and if $\theta>\aleph_0$ we also require that this sequence is
continuous. In addition, we require that
$\xi_{\epsilon+1}$ is a successor ordinal,
for every $\epsilon<\theta$.

 For every $\alpha\in\kappa^+$, we define
 $$c_\alpha\deq\{\epsilon<\theta:\,(\xi_\epsilon,\xi
 _{\epsilon+1})\cap
 E_\alpha\neq\emptyset\},$$
 $$d_\alpha\deq\{\sup(\xi_{\epsilon+1}\cap E_\alpha):\,\epsilon\in
 c_\alpha\}.$$

Note that $d_\alpha\subseteq E_\alpha
\setminus\{\xi_\epsilon:\,\epsilon<\theta\}$ and $\sup(d_\alpha)=\kappa$,
for each $\alpha<\kappa^+$.

 Now we can distinguish two cases:

 $\underline{\hbox{Case (A)}.}$ For some $\alpha$, for every $\beta\in
 (\alpha,\kappa^+),$ the symmetric difference $d_\beta\symmdif d_\alpha$ is 
bounded
 in $\kappa$.

 We set $C^{\ast}=d_\alpha$ for some such $\alpha$, and easily check that
 (1) holds.

 $\underline{\hbox{Case (B)}.}$ Not (A).

 Therefore, for every $\alpha\in\kappa^+$, there is a minimal
 $f(\alpha)\in
 (\alpha,\kappa^+)$ such that $\kappa=\sup(d_\alpha\symmdif
 d_{f(\alpha)}).$
The following can be noticed:

{\it Observation.\/}
 $$f(\alpha)\leq\beta<\kappa\implies
 \kappa=\sup(d_\alpha\symmdif
 d_\beta).\eqno(\ast)$$
[Why?
If $\gamma$ is such that $(E_{f(\alpha)}\setminus E_\alpha)
\cup (E_\beta\setminus E_{f(\alpha)})\subseteq\gamma$, and
$\epsilon\in c_\beta\cap c_\alpha$ is such that
$\sup(E_\beta\cap\xi_{\epsilon+1})=\sup(E_\alpha\cap\xi_{\epsilon+1})$
and $\xi_\epsilon\ge\gamma$, then
we also have $\sup(E_\alpha\cap\xi_{\epsilon+1})=\sup(E_{f(\alpha)}\cap
\xi_{\epsilon+1}).$]

 Let
 $$E^{\ast}\deq\{\delta<\kappa^+:\,\delta\hbox{ is a limit ordinal such that }
 \bwedge_{\alpha<\delta}f(\alpha)<\delta\}.$$
 Obviously, $E^{\ast}\in V_2$ is a club of $\kappa^+$.
 Then by $(\delta)$, for some $\alpha^\ast\in S_{i(\ast)}$ we have 
$$V_2\satisfies\,``\cf(\alpha^\ast)=\mu^\ast\,\,\&\,\,C_{\alpha^\ast}^{i(\ast)}
\cap E^{\ast}\hbox{ is a club of }\alpha^\ast".$$

 Let $\langle\beta_i:\,i<\rho\rangle$ be an increasing enumeration of
${\rm acc}(C_{\alpha^\ast}^{i(\ast)}\cap E^{\ast})$,
for some $\rho\ge\mu^{\ast}$.
We construct by induction on $i$
an increasing subsequence $\langle\xi_{\epsilon_i}:\,
i<\mu^{\ast}\rangle$ of $\langle\xi_\epsilon:\,
\epsilon<\theta\rangle$, thus obtaining a contradiction
with $\mu^{\ast}>\theta.$

Let us introduce the notation $\delta^\alpha_\epsilon\deq
\sup(\xi_{\epsilon+1}\cap E_\alpha)$ for $\epsilon\in c_\alpha$
and $\alpha\in\kappa^+$.

We first make two observations:

\item{$(o_1)$} If $i<j<\rho$, then by
the square property of $C_{\alpha^\ast}^{
i(\ast)}$ and (iii) above, $E_{\beta_i}\supseteq E_{\beta_j}$. So, if
$i<j<k<\rho$ and $\epsilon<\theta$ is such that
$(d_{\beta_i}\symmdif d_{\beta_j})\cap(\xi_{\epsilon},\xi_{\epsilon+1})
\neq\emptyset$ and $(d_{\beta_j}\symmdif d_{\beta_k})\cap
(\xi_{\epsilon},\xi_{\epsilon+1})\neq\emptyset$, then
$\delta^{\beta_i}_\epsilon>\delta^{\beta_j}_\epsilon>\delta^{\beta_k}_\epsilon.$

\item{$(o_2)$} Given $i,j<\rho$ (so 
$f(\beta_i)< \beta_j)$, and $\epsilon<\theta$, there is a
$\zeta>\epsilon$ such that $(d_{\beta_i}\symmdif d_{\beta_j})
\cap (\xi_\zeta,\xi_{\zeta+1})\neq\emptyset$. This follows from the
$(\ast)$ above.

Now we proceed 
to
construct $\langle\xi_{\epsilon_i}:\,i<\mu^\ast\rangle$.

Let $\gamma_0^0=\beta_0$ and $\gamma_1^0=\beta_1$
(so $\gamma^0_1>
f(\gamma_0^0)$) and let $\epsilon_0$ be such that
$(d_{\gamma^0_0}\symmdif d_{\gamma_1^0})\cap (\xi_{\epsilon_0},
\xi_{\epsilon_0+1})\neq\emptyset$.

Suppose that for some $i<\mu^\ast$, we have chosen $\gamma^i_0,
\gamma^i_1$ and $\epsilon_i$ so that $(d_{\gamma_0^i}
\symmdif d_{\gamma_1^i})\cap (\xi_{\epsilon_i},
\xi_{\epsilon_i+1})\neq\emptyset.$ We 
wish to define $\epsilon_{i+1}$. Using the chosen
$\gamma^i_0,\gamma^i_1$, we build an increasing
subsequence $\gamma^i_0,\gamma_1^i,
\gamma_2^i\ldots$ of $\langle\beta_k:\,k<\rho\rangle$
such that
$$l_1<l_2\implies (\xi_{\epsilon_i},\xi_{\epsilon_i+1})
\cap (d_{\gamma^i_{l_1}}\symmdif d_{\gamma^i_{l_2}})\neq\emptyset.$$
By $(o_1)$ above, the sequence $\gamma^i_0,\gamma^i_1,
\gamma^i_2\ldots$
must stop after a finite number of stages, since
otherwise we obtain an
infinite decreasing sequence of ordinals. Let $\gamma^{i+1}_1>
\gamma^{i+1}_0>\max\{\gamma^i_l:\,\gamma^i_l\hbox{ is well defined}\}$ be
such that $\gamma_0^{i+1}=
\beta_{k_0}$ and $\gamma_1^{i+1}=\beta_{k_1}$ for some $k_0<k_1<\rho.$
Let $\epsilon_{i+1}
>\epsilon_i$ be
such that $(d_{\gamma^{i+1}_0}\setminus d_{\gamma^{i+1}_1})
\cap(\xi_{\epsilon_{i+1}},\xi_{\epsilon_{i+1}+1})
\neq\emptyset.$ Such an $\epsilon_{i+1}$ exists by $(o_2)$ above.

If $i$ is a limit ordinal $<\mu^{\ast}$, we define $\gamma^i_0>
\sup\{\gamma^j_l:\,j<i\,\,\&\,\,\gamma^j_l\hbox{ is well defined}\}$,
$\gamma_0^{i}=\beta_k$ for some $k<\rho$, which is possible
by
the regularity of $\mu^\ast>\aleph_0$.
We define $\gamma^{i}_1$ and $\epsilon_i$
as above.

 \smallskip

 (2) Let $\gamma<\kappa$ be fixed. In $V_2$,
 let $\sigma$ be regular, $\sigma\in(\gamma,
 \kappa)$, so $\sigma$ is also regular in $V_1$. We follow the ideas of
 [Sh 365,\S2] but the following is self-sufficient.

 The plan is to replace in the proof of (1), each $d_\alpha$
 by a somewhat larger set $d_\alpha^{\dagger}$, so that
 $$\beta\in\nacc(d_\alpha^{\dagger})\implies{\cf}^{V_2}(\beta)\ge\sigma.$$
 In the proof of (1), $C^\ast=d_\beta$ for some $\beta$. The
 definition
 of $d_\alpha^{\dagger}$'s will be such that putting $C^\ast=d_\beta^{\dagger}$
 for the
 same
 $\beta$, the newly defined $C^{\ast}$ will still satisfy (1).

 First, we shall
 enlarge each $d_\alpha$ to a $d_\alpha^{\dagger}$, so that
 $$\beta\in\nacc(d_\alpha^{\dagger})\,\&\,{\cf}^{V_2}(\beta)<\sigma\implies\beta\in
 \nacc(E_\alpha).
 \eqno(\ast\ast)$$
 In $V_1$, we choose for each $\delta\in\kappa$ a club $e_\delta$ in $\delta$
 such that
 $\otp(e_\delta)=\cf^{V_1}(\delta)$, and define $\bar{e}\deq\langle e_\delta:\,
 \delta<\kappa\rangle\in V_1$.

 Fix an $\alpha<\kappa^+$ and set $C\deq\{\xi_{\epsilon+1}:\,\epsilon\in c_\alpha\}$,
 $E\deq E_\alpha$, $d\deq d_\alpha$.
 Then, in the notation of [Sh 365,\S2],
using the ``glue'' operator $\gl$,
 $$d={\gl} (C,E)\deq {\gl}^0 (C,E)\deq {\gl}^1_{\sigma,0} (C,E,\bar{e})\deq
 \{\sup(\beta\cap E):\,\beta\in C,\,\beta>\min(E)\}.$$
 We shall have
 $$d_\alpha^{\dagger}={\gl}^1_\sigma (C,E,\bar{e})\deq\cup_{n\in\omega}
 {\gl}^1_{\sigma,n}(C,E,\bar{e}),$$
 where ${\gl}^1_{\sigma,n}$ are defined inductively on $n$ as follows
 (${\gl}^1_{\sigma,0} (C,E,\bar{e})$ is given above):

 $${\gl}^1_{\sigma,n+1} (C,E,\bar{e})\deq{\gl}^1_{\sigma,n} (C,E,\bar{e})\cup$$
 $$\eqalign{\cup\Bigl\{\sup(\delta\cap E):\hbox{ for some
 }\beta\in\nacc\bigl({\gl}^1_{\sigma,n} \bigl(C,E,\bar{e})\bigr)
\hbox{ we have:}\cr
{\cf}^{V_2}(\beta)<
\sigma\,\&\,\delta\in e_\beta
\,\&\,
 \delta>\sup\bigl(E\cap\delta\cap{\gl}^1_{\sigma,n} (C,E,\bar{e})\bigr)
 \Bigr\}.\cr}$$
 We can easily check that this definition leaves us in the
 situation of $(\ast\ast)$.

 On the other hand,
since
 $E\in V_1$ is a club of $\kappa>\sigma$, then
 $$E^{\dagger}\deq\{\delta\in E:\,\otp(\delta\cap E)\hbox{ is divisible by
 }\sigma\}$$
 is a club of $\kappa$, is in $V_1$ (because
$E\in V_1$), and has the property that
 $$\beta\in\nacc(E^{\dagger})\implies{\cf}^{V_2}(\beta)\ge\sigma.$$
(Of course, also ${\cf}^{V_1}(\beta)\ge\sigma.$)

 Now,
looking at all $\alpha<\kappa^+$ simultaneously,
had we in our definition of $\langle
 E_\alpha:\,\alpha<\kappa^+\rangle$ in
 the proof of (1)
 replaced each $E_\alpha$ by $E_\alpha^{\dagger}$, (i)-(iii)
would still have been
 satisfied.
So, by $(\ast\ast)$, had we used
$E_\alpha^{\dagger}$ rather than $E_\alpha$ in the definition of
 $d_\alpha^{\dagger}$,
we would have have:
 $$\beta\in\nacc(d_\alpha^\dagger)
\implies{\cf}^{V_2}
(\beta)\ge\sigma,$$
as desired.

 It only remains to check that $C^\ast$
would still satisfy (1),
which is easy.

(3) Similar to (2).\eop${}_{2.0.}$

\medskip

{\bf Remark 2.1}(0) Note that a particular case of {\it (iv)\/}
of Theorem 2.0
is that $V_1\satisfies\,2^\kappa=\kappa^+$.

(1) In fact, {\it (iv)\/} could be weakened using [Sh 430, \S1].

\medskip

Now we shall see that the same conclusion as in Theorem 2.0(1) can be
obtained under somewhat different assumptions. We have
no assumptions on $\cf(\Club^{V_1}(\kappa),\supseteq)$
or $\kappa^+$, but we have to
add some assumptions on certain covering numbers. If we use
forcing to go from $V_1$ to $V_2$, these correspond to the chain
condition of the forcing using. We also do not need to require that
$\kappa$ is a limit cardinal of $V_1$.

\medskip

 {\bf Theorem 2.2.} Assume that:

\item{{\it (i)\/}
} $V_1\subseteq V_2$ are models of $ZFC$, $V_1$ is a class of $V_2$
containing all ordinals
and $V_1,V_2$ are transitive.

\item{{\it (ii)\/}
} $V_1\satisfies``\kappa\hbox{ is regular }>\aleph_0".$

\item{{\it (iii)\/}
} $V_2\satisfies``\cf(\kappa)=\theta,\,\kappa\hbox{ is a cardinal and }
\theta<\kappa".$

Moreover, assume that:

$\underline{\hbox{either}}$ $\theta>\aleph_0$ and

\item{{\it (iv${}^{'}$)\/}} {\it If} $a\in V_2$ is such that 
$$V_2\satisfies``\card{a}
\le\cf(\Club^{V_2}(\theta),\supseteq)",$$
and 
$$V_2\satisfies
``a\subseteq\Club^{V_1}(\kappa)",$$

{\it then\/} we can find a sequence $\langle b_i:\,i<\theta\rangle\in V_2$
such that $\bwedge_{i<\theta} b_i\in V_1$,
while \break $V_1\satisfies``\card{b_i}\le\kappa"$
and $V_2\satisfies``a\subseteq\cup_{i<\theta} b_i".$
(Hence, $V_2\satisfies ``{\rm cf}({\rm Club}^{V_2}(\theta),\supseteq)
\le\kappa"$.)

$\underline{\hbox{or}}$ $\theta=\aleph_0$ and:

\item{{\it (iv${}^{''}$)\/}} $V_2\satisfies\,``2^{\aleph_0}<\kappa".$
\item{{\it (v${}^{''}$)\/}} For every $a\subseteq
\Club^{V_1}(\kappa)$ such that $\card{a}^{V_2}
\leq 2^\theta
$, there is a $b\in V_1$ such that $a\subseteq b$ and
${\card{b}}^{V_1}\le\kappa$.

{\it Then\/}
in $V_2$, we can find an
unbounded $C^\ast\subseteq \kappa$
 such that
 $$E\in \Club^{V_1}(\kappa)\implies C^\ast\setminus E\hbox{
 is bounded.}$$



\medskip

Before proving this theorem, let us make some remarks.

\medskip

{\bf Remark 2.3}(0) Suppose that
the first three assumptions of 2.2
are satisfied.

If ${\rm cf}({\rm Club}^{V_2}(\theta),
\supseteq)<\kappa$ and $V_2$ is obtained from $V_1$ by
changing the cofinality of $\kappa$ to $\theta$ via a
$\kappa^+-cc$ forcing,
{\it then (iv${}^{'}$)\/} holds (in either of the two cases of the theorem).
Actually, {\it (iv${}^{'}$)\/} holds under some weaker conditions
(see [Sh 410, \S2] for this). 

(1) It is also meaningful to use
the notion of $\cf(\Club(\theta),\supseteq)$
for $\theta=\aleph_0$. Namely a club subset of $\omega$ is simply any unbounded
set, and then the $\cf(\Club(\aleph_0),\supseteq)$ corresponds to the familiar
notion of ${\bf d}$, the smallest cardinality of a dominating family in
$({}^\omega\omega,\le^{\ast}).$ If {\bf d}$=2^{\aleph_0}$, then obviously
{\it (iv${}^{'}$)\/} is weaker than {\it (v${}^{''}$).\/}

(2) If $V_1\satisfies``\cf(\Club(\kappa),\supseteq)=\kappa^+"$ (as in Theorem 2.0),
then {\it (v${}^{''}$)\/} holds. The natural case
that we have in mind is $V_2\satisfies\, 2^\theta<\kappa$.
So, if we are using a $\kappa^+$-cc forcing to pass from $V_1$ from $V_2$
and  $V_2\satisfies 2^\theta<\kappa$,
theorem 2.2 is stronger than 2.0(1).

(3) Note that if $V_1\satisfies``\kappa=\mu^+$ and $\mu$ is regular"
then $V_2\satisfies``\cf(\kappa)=\cf(\mu)"$, by [Sh -g VII \S4]=
[Sh -b XIII \S4].

(4) We mentioned earlier, a different proof of 2.2
can be found in [Gi1], with somewhat stronger assumptions in the case of $\theta>\aleph_0$
(that is, assuming $2^\theta <\kappa)$.
\medskip

We proceed to the proof of 2.2.

\medskip

{\bf Proof of 2.2.} We break the proof into two cases:

$\underline{\hbox{ Case 1.}}$ $\theta>\aleph_0$.

Let $d=\{\alpha_\zeta:\,\zeta<\theta\}\in V_2$ be a club of $\kappa$.
For each club $A\subseteq\theta$ from $V_2$, we try to set
$C^{\ast}=d_A\deq\{\alpha_\zeta:\zeta\in A\}$
and obtain $C^{\ast}$ which satisfies the theorem.
If we succeed for some $A$, the theorem is
proved. Note that each $d_A$ is unbounded in $\kappa$.

Otherwise, for each $A$ as above,
we can choose an $E_A\in\Club^{V_1}(\kappa)$,
which witnesses that $C^{\ast}=d_A$ does not work. So, in $V_2$
we have
$\theta=\sup\{\zeta\in A:\,\alpha_\zeta\notin E_A\}.$
Without loss of generality, each $E_A$ is a club subset of $\kappa$ in $V_1$.

Let $\PP\subseteq\Club^{V_2}(\theta)$ be cofinal
in $({\rm Club}^{V_2}(\theta),\supseteq)$, of cardinality
$\cf^{V_2}(\Club(\theta))$.
Therefore $a\deq\{E_A:\,A\in\PP\}$ is a subset of $\Club^{V_1}(\kappa)$
of cardinality $\le\cf^{V_2}(\Club(\theta),\supseteq)$. So we
can find $\langle b_i:\,i<\theta\rangle$ as guaranteed by {\it (iv${}^{'}$)\/}.
 
In $V_1$, for $i<\theta$, let $\{C:\, C\in b_i\hbox{ and }
C\hbox{ is a club of }\kappa\}$ be enumerated as $\{C_j^i:\,j<j_i\le\kappa\}$
and let
$E_i\deq
\Delta_{j<j_i}C_j^i$. In $V_2$, the set $E_i\cap d$ is a club of $\kappa$.
Therefore, $A_i\deq\{\zeta<\theta: \alpha_\zeta\in E_i\}$ is a club of
$\theta$.

Now consider $A^{\ast}\deq\Delta_{i<\theta} A_i=\{\zeta<\theta:\,
\bwedge_{i<\zeta}\alpha_\zeta\in E_i\}$. It is a club of $\theta$, so we can
find an $A^\tensor\in\PP$ such that $A^\tensor\subseteq A^\ast$.
Without loss of generality, $A^\tensor$ is a club.
So, $E_{A^\tensor}$ is well defined, and for some $i<\theta$, $E_{A^\tensor}
\in b_i$, so $E_{A^\tensor}
=C_j^i$ for some $j<\kappa$. Fix such $i$ and $j$.

Then by the definition of $E_i$ and $A^\ast$, we have
$$\zeta\in A^{\ast}\,\&\,
\zeta>i\,\,\&\,\,\alpha_\zeta>j\implies \alpha_\zeta\in E_{A_\tensor},$$
which is a contradiction
with $A^\tensor\subseteq A^{\ast}$ and $\theta=\sup\{\zeta\in A^\tensor:\,
\alpha_\zeta\notin E_{A^\tensor}\}$.

\smallskip

$\underline{\hbox{Case 2.}}$ $\theta=\aleph_0.$

As before, we first fix (in $V_2$) an increasing unbounded
sequence $d=\{\alpha_n:\,n\in\omega\}$ in $\kappa$.
For all $E\in\Club^{V_1}(\kappa)$, we try
setting
$$C^\ast=C^\ast(E)=\bigl\{\sup(\alpha_n\cap E):\,n\in\omega\bigr\}.$$
Note that each $C^\ast(E)$ is unbounded in $\kappa$.
If we succeed for some $E$, then we are done.
So, from now on we assume otherwise.

To
do the proof, we need the following fact.

\smallskip

{\bf Fact 2.2.a.} If $\langle E_i:\,i<i^\ast<(2^{\aleph_0})^+\rangle\in V_2$
is such that each
 $E_i\in\Club^{V_1}(\kappa)$,
{\it then\/} we can find an $E\in\Club^{V_1}(\kappa)$ such that
$$\cwedge_{i<i^\ast}\sup(E\setminus E_i)<\kappa.$$

{\bf Proof.} Without loss of generality, each $E_i$ is a club of $\kappa$
in $V_1$. Now 
consider $a\deq\{E_i:\,i<i^\ast\}$, and apply
{\it (v${}^{''}$)\/} to find a $b\in V_1$ such that ${\card{b}}^{V_1}\leq\kappa$
and $a\subseteq b$. In $V_1$, by intersecting $b$ with
the family of all clubs of $\kappa$ in $V_1$
if necessary, we can without loss of generality assume that
each element of $b$ is a club of $\kappa$ in $V_1$.
So, let us in $V_1$ enumerate $b$ as $b=\{C_\epsilon:\,
\epsilon<\kappa\}$. Then let $E\deq\Delta_{\epsilon<\kappa}C_\epsilon$, so $E\in\Club^{V_1}(\kappa).$
It is easily seen that
$E$ is as required.\eop${}_{2.2.a.}$

\smallskip

{\it Continuation of the proof of 2.2, Case 2.\/}
We place ourselves in $V_2$ and choose by induction on
$i<(2^{\aleph_0})^+$, an $E_i\in\Club^{V_1}(\kappa)$ such that:

\item{$(\alpha)$} $j<i\implies E_j\supseteq^{\ast} E_i$.

\item{$(\beta)$} $E_{i+1}\subseteq E_i$ shows that setting
$C^\ast=C^\ast(E_i)$
does not give a satisfactory definition of $C^\ast$.

This induction is easily carried: at the successor stages we do
as in $(\beta)$, and at the limit stages we
use Fact 2.2.a.

For a given $i<(2^{\aleph_0})^+$, if $n\in\omega$ is such that
$\sup(\alpha_n\cap E_i)\notin E_{i+1}$, then
$\sup(\alpha_n\cap E_i)>\sup(\alpha_n\cap E_{i+1}).$
Then it follows from $(\beta)$ that
$$\cwedge_{i<(2^{\aleph_0})^+}
\{n\in\omega:\,\sup(\alpha_n\cap E_{i+1})<\sup(\alpha_n\cap E_i)\}
\hbox{ is infinite}.$$
Now we define for each $i<(2^{\aleph_0})^+$, a function $h_i:\omega\to\kappa$ by
$$h_i(n)=\sup(\alpha_n\cap E_i).$$
So we have:
\item{$(\alpha)^\ast$} $j<i
<(2^{\aleph_0})^+\implies h_j \le_{J^{\bd}_\omega} h_i$.

\item{$(\beta)^\ast$} $h_{i+1}\neq h_i (\mod J^{\bd}_\omega)$. 

This is a contradiction.
\eop${}_{2.2.}$

\medskip

Similar in its proof to Theorem 2.2, the following theorem enables us to
do outside guessing of clubs using a proper filter in $V_2$.

\medskip

 {\bf Theorem 2.4.} Assume that:

 \item{{\it (i)\/}} $V_1\subseteq V_2$ are transitive classes containing the ordinals
and satisfying $ZFC$.

 \item{{\it (ii)\/}} $V_1\satisfies
``\kappa\hbox{ is a regular cardinal }>\aleph_0".$

 \item{{\it (iii)\/}} $V_2\satisfies
``\kappa\hbox{ is a singular cardinal},\cf(\kappa)=
\theta"$.

 \item{{\it (iv)\/}} $V_2\satisfies\,``\theta>\aleph_0".$

 \item{{\it (v)}}
 $V_1\satisfies\,``\mu=\cf\bigl(\Club(\kappa),\supseteq\bigr)"$, and

\item{{\it (vi)\/}} $\sigma$ is a cardinal in $V_2$ and
 $$\Bigl(\forall a\Bigr)\Bigl[a\in V_2\,\&\,
 a\subseteq\mu\,\&\,
 \card{a}<\sigma\to\bigl(\exists b\in V_1\bigr)\bigl(a\subseteq
 b\subseteq\mu\,\&\,\card{b}
 <\kappa\bigr)\Bigr].$$

 {\it Then\/} we can find in $V_2$ an increasing continuous sequence
 $\langle
 \alpha_\zeta^{\ast}:\zeta<\theta\rangle$ with limit $\kappa$, and a
 proper $\sigma$-complete filter $\DD$ on $\theta$, such that:
 $$\hbox{For every club }E\hbox{ of }\kappa\hbox{ from }V_1,
\hbox{ we have }
 \{\zeta:\,\alpha^\ast_{\zeta+1}\in E\}\in\DD.\eqno{(\ast\ast\ast)}$$

\medskip

{\bf Remark 2.5}(0) If $\sigma=\aleph_0$, {\it (vi)} is redundant.

(1) The assumptions imply that $\sigma\leq\theta$. Otherwise we could use
{\it (vi)\/} on a cofinal $\theta$-sequence of $\kappa$ in $V_2$, and obtain a
contradiction with $V_1\satisfies``\kappa$ is regular''.

(2) If instead of {\it (vi)\/}
we have other properties of the style of [Sh 355,\S5],
or {\it (iv${}^{'}$)\/} of Theorem 2.2, then we get corresponding completeness properties
of the filter. For example, we could have that among any $\rho$ members of
$\DD$, there are $\lambda$ whose intersection is also in $\DD$, for some
cardinal $\rho\ge\lambda$.

(3) Remark 2.3.(3) applies here too.

\medskip

 {\bf Proof of 2.4.} Once we define $\langle
 \alpha_\zeta^{\ast}:\zeta<\theta\rangle$, we shall have that
 $$\DD=\Bigl\{A\subseteq\theta:\,\hbox{for some
 }E\in\Club^{V_1}(\kappa),\hbox{ we have that }
 \bigl(\alpha^\ast_{\zeta+1}\in E\implies\zeta\in A\bigr)\Bigr\}.$$
 This definition makes sense for any $\langle
 \alpha_\zeta^{\ast}:\zeta<\theta\rangle$ increasing continuous
 with limit $\kappa$ and yields a $\sigma$-complete
 filter $\DD$ (by {\it (vi)\/}),
 and $(\ast\ast\ast)$ holds. The point is to have that $\DD$ is
 proper,
 i.e. $\emptyset\notin\DD$, and we now show how to make the
 choice
 of $\langle
 \alpha_\zeta^{\ast}:\zeta<\theta\rangle$ which will satisfy this.

In $V_2$ we
 let $\{\xi_\epsilon:\,\epsilon<\theta\}$ be a strictly increasing cofinal
sequence in $\kappa$.
 For each club $E$ of $\kappa$ from $V_1$, let

 $$d_E\deq\{\sup(\xi_\epsilon\cap E):\,\epsilon<\theta\,\&\,
 E\cap(\xi_\epsilon,
 \xi_{\epsilon+1})\neq\emptyset\,\&\,\xi_\epsilon
>\min(E)\}.$$
So, every $d_E$ is an unbounded subset of $\kappa$.
 If for some $E$, an increasing enumeration of $d_E$ can be used for
 $\langle
 \alpha_\zeta^{\ast}:\zeta<\theta\rangle$, and a proper
 filter is obtained, we are done.

 Otherwise, let $\{E_\alpha:\,\alpha<\mu\}$ be a cofinal
sequence of clubs of $\kappa$ in $V_1$
and assume that the enumeration is in $V_1$. By induction on $n\in\omega
$, we
 choose
 $\alpha_n<\mu$ such that:

 \item{-}$\alpha_0=0.$

\item{-}$\alpha_{n+1}$ is such that $E_{\alpha_{n+1}}$ is a club
 exemplifying the
 failure of $d_{E_{\alpha_n}}$ to give a satisfactory definition of
 $\langle
 \alpha_\zeta^{\ast}:\zeta<\theta\rangle$. Without loss of generality,
$E_{\alpha_{n+1}}\subseteq E_{\alpha_n}$.

 Then $E^{\ast}=\cap_{n\in\omega}E_{\alpha_n}$ is an element of $V_2$,
 and a club
 of $\kappa$, so necessarily
$$\theta=\sup\{\epsilon:\,(\xi_\epsilon,\xi_{\epsilon+1})
 \cap E^{\ast}\neq\emptyset\}.$$
Note that the definition of $E_{\alpha_{n+1}}$ implies the existence of an
$\epsilon_n(\ast)<\theta$ such that
 $$\epsilon\in(\epsilon_n(\ast),\theta)\,\&\,
 (\xi_\epsilon,\xi_{\epsilon+1})\cap E_{\alpha_{n+1}}\neq\emptyset
 \implies
 \sup\bigl((\xi_\epsilon, \xi_{\epsilon+1})\cap E_{\alpha_{n+1}}\bigr)
 <\sup\bigl((\xi_{\epsilon},\xi_{\epsilon+1})\cap E_{\alpha_n}\bigr).$$
 
 Now  $\cup_{n\in\omega}\epsilon_n(\ast)<\theta$,  so take an
 $\epsilon\in(\cup_{n\in\omega}\epsilon_n(\ast),\theta)$
 such  that
 $((\xi_\epsilon,\xi_{\epsilon+1})\cap E^{\ast}\neq\emptyset.$
 Then the sequence
 $\bigl\langle\sup\bigl(\xi_\epsilon,\xi_{\epsilon+1})
 \cap E_{\alpha_n}\bigr):\,n\in\omega\bigr\rangle$ is a strictly
 decreasing sequence
 of ordinals, a contradiction.\eop${}_{2.4}$

\bigskip

{\bf \S3 The family} $I_{<f}[\lambda]${\bf .} Throughout this paper, we dealt
with sets with squares, some weak version of squares, and the ideal $I[\lambda]$.
In this section we would like to explore the connection between these notions.
We show that the notion of square, the one of the weak square, and the
ideal $I[\lambda]$ are all definable by the same definition, simply
by changing a certain parameter. We shall also discuss general properties
of notions which can be defined in this way.

We now review some familiar notions and related families of sets. 

To avoid trivialities, in the following $
\lambda$ is always assumed to be
a regular uncountable cardinal.


\smallskip

{\bf Definition 3.0}(0) For a subset $S^+$ of $\lambda$, we say that
$S^+$ {\it has a square\/}, if the following holds:

There is
a sequence $\langle C_\delta:\,\delta\in S^+\rangle$
such that:

\item{$(a)$} $C_\delta$ is a closed subset of $\delta$. 

\item{$(b)$} $\beta\in C_\delta\implies 
\beta\in S^+\,\,\&\,\, C_\beta=\beta\cap C_\delta$.

\item{$(c)$} $\beta$ is a limit ordinal $\in S^+\iff \beta=
\sup(C_\beta).$

\item{$(d)$} $\otp(C_\delta)<\delta$.

$I^\square[\lambda]$ is the family of all subsets $S$ of $\lambda$
for which there is an $S^+\subseteq\lambda$ which has a square
and satisfies $S\setminus S^+$ is non-stationary.

(1) A subset $S$ of $\lambda$ is said to 
{\it have a weak square\/} if $S\subseteq S^+ \bigl(\mod
\Club(\lambda)\bigr)$ for some $S^+\subseteq\lambda$ with the
following property:

There is a sequence $\langle \PP_\delta:\,\delta\in S^+\rangle$
such that:

\item{$(a)$} Each $\PP_\delta$ is a family of closed subsets of $\delta$,
and if $\delta$ is a limit ordinal, all members of $\PP_\delta$ are unbounded
in $\delta$.

\item{$(b)$} $\beta\in a\in\PP_\delta\implies a\cap\beta\in\PP_\beta.$

\item{$(c)$} $\card{\PP_\delta}\leq\card{\delta}.$

\item{$(d)$} $a\in\PP_\delta\implies\otp(a)<\delta.$

$I^{\rm w\square}[\lambda]$ is the family of all subsets of $\lambda$
which have a weak square.

(2) For an
$S\subseteq\lambda$, we say that
$S$ is {\it good\/}, if $S\subseteq S^+\bigl(\mod \Club(\lambda)\bigr)$
for some $S^+\subseteq\lambda$
such that there
is a sequence
$\bar{C}=\langle C_\delta:\,\delta<\lambda\rangle$ for which: 

\item{$(a)$} $C_\delta$ is a closed subset of $\delta$.

\item{$(b)$} If $\alpha\in\nacc(C_\delta)$, then $C_\alpha=
\alpha\cap C_\delta$.

\item{$(c)$} If $\delta\in S^+$, then $\delta=\sup(C_\delta)$
and $\otp(C_\delta)=\cf(\delta)<\delta.$

\item{$(d)$} $\nacc(C_\delta)$ is a set of successor ordinals.


$I[\lambda]$ is the family of good subsets of $\lambda$.

\smallskip

{\bf Remark 3.1}(0) The notions of square and weak square are well
known and were introduced by Jensen.
The first appearance of $I[\lambda]$ is in
[Sh 108], or [Sh 88a]. 
It was also considered in [Sh 351], [Sh 420] and elsewhere.
We have already used
$I[\lambda]$ in the first section. 
The definition we use differs from the original definition from [Sh 88a]
for example, but the equivalence is proved in [Sh 420, 1.2]. It
is shown in [Sh 88a, 3(1)] that $I[\lambda]$ is a normal ideal on $\lambda$.
Under certain circumstances, like when $\lambda$
is the successor of a strong limit singular, the ideal $I[\lambda]$ contains
a maximal set [Sh 108], and we have made use of this fact
in the first section. 
For various further properties of $I[\lambda]$ see the above references.

(1) Obviously, $I[\lambda]\cap I^\square[\lambda] \cap I^{{\rm w}\square}
[\lambda]\supseteq NS[\lambda]$, and each $I[\lambda], I^{{\rm w}\square}
[\lambda]$ and $I^\square[\lambda]$ is closed under taking subsets. 

(2) The notation $I^\square[\lambda]$ and $I^{{\rm w}\square}[\lambda]$ should
not suggest that these families are necessarily ideals.

(3) Note that what we have defined now as a square on $S$
differs from the definition we used in $\S1$. This does not have any effect
on $I_1[\lambda]$ (see below),
so we adopt the present definition for convenience.

\smallskip

We now introduce a notion which as
its particular cases includes $I[\lambda]$, $I^\square[\lambda]$
and
\break $I^{{\rm w}\square}[\lambda]$.

\smallskip

{\bf Definition 3.2}(0) Let $\lambda$ be a regular uncountable
cardinal, and
$f$ a function from $\lambda$ to the cardinals.
To avoid trivialities, we assume that $f(i)\ge 2$, for all $i\in\lambda$.

We define
$$I_{<f}[\lambda]$$
as the family of all $S\subseteq\lambda$ for
which
there is an $S^+\subseteq\lambda$ such that
$S\setminus S^+$ is non-stationary, and a
sequence
$\bar{C}=\langle C_\delta:\,
\delta\in S^+\rangle$ which has the following properties:

\item{$(a)$} $C_\delta$ is a closed subset of $\delta$.
\item{$(b)$} $\otp(C_\delta)<\delta$.
\item{$(c)$} $\delta\in S^+$ is a limit ordinal $\implies\delta
              =\sup(C_\delta)$.
\item{$(d)$} $\nacc(C_\delta)
$ is a set of successor ordinals.
\item{$(e)$} For every $\beta<\lambda$, the set
$$\PP_\beta=\PP_\beta[\bar{C}]\deq\{C_\delta\cap\beta:\,\beta\in C_\delta\}$$
has cardinality $<f(\beta)$.

(1) We call $\bar{C}$ and $S^+$ as above
{\it witnesses} for $S\in I_{<f}[\lambda]$. Or, we say that they
{\it exemplify\/} that $S\in I_{<f}[\lambda].$

(2) If $f(i)=g(i)^+ $ for all $i$, we let
$$I_g[\lambda]=I_{<f}[\lambda].$$
If $f$ is constantly equal to $\mu$, we write $I_{<\mu}[\lambda]$\
(or $I_\mu[\lambda]$) instead of $I_{<f}[\lambda]$ (or $I_f[\lambda]$).

\smallskip

{\bf Remark 3.3}(0) It is easy to see that for any choice of $f$,
the family
$I_{<f}[\lambda]$ includes all non-stationary subsets of $\lambda$,
and that $I_{<f}[\lambda]$ is closed under taking subsets.

(1) Again, $I_{<f}[\lambda]$ is not always an ideal, and we shall below
discuss sufficient conditions for $I_{<f}[\lambda]$ to
be an ideal, or a normal ideal. 

(2) By $(b)-(c)$ of the Definition 3.2.0, if $S\in I_{<f}[\lambda]$, then
$S\cap REG$ is not stationary.

(3) Notice that
if $f$ and $g$ are functions from $\lambda$ into the cardinals $\ge 2$,
and $g(i)\le f(i)$ for all $i$, then $I_{<g}[\lambda]\subseteq
I_{<f}[\lambda].$

\smallskip

We first make some general remarks about $I_{<f}[\lambda]$ which show to which
extent we can generalize results about $I[\lambda]$ to this notion. 
The following notation will be convenient,
and corresponds to the one introduced in [Sh 420]
to study $I[\lambda]$.

\smallskip

{\bf Notation 3.4.}  Suppose that $\lambda$ is a cardinal, $S\subseteq\lambda$
and $f$ is a function from $\lambda$ into the cardinals.
The symbol
\def\ooplus{{}^-\oplus^f_{\bar{\PP},S^+}(S)}
$\ooplus$ means that $S^+\subseteq\lambda$
and $\bar{\PP}=\langle\PP_\alpha:\,\alpha\in\lambda\rangle$ is a sequence
of sets with the following properties:

\item{$(i)$} $\PP_\alpha$ consists of $<f(\alpha)$ subsets of $\alpha$.

\item{$(ii)$} If $\delta\in S^+$ is a limit ordinal, then there is an $x
\subseteq\delta$ such that
$$\otp(x)<\delta=\sup(x)\quad\&\quad\bwedge_{\zeta<\delta}x\cap\zeta\in
\bigcup_{\gamma<\delta}\PP_\gamma.$$

\item{$(iii)$} $S\setminus S^+$ is non-stationary.

\smallskip

We formulate
a Lemma which corresponds to the fact
[Sh 420, \S1],
that the various
definitions of $I[\lambda]$ considered in [Sh 88a], [Sh 108] and [Sh 420]
are equivalent.

\medskip

{\bf Lemma 3.5.} Suppose that $\lambda$ is an uncountable regular
cardinal and $f$ is a function from $\lambda$ into the cardinals.
Let us enumerate the following statements:

\item{(1)} $S\in I_{<f}[\lambda]$.

\item{(2)} $\ooplus$ for some $\bar{\PP}$ and $S^+$ such that all elements of
$\bigcup_{\alpha<\lambda}\PP_\alpha$ are closed.

\item{(3)} $\ooplus$ for some $\bar{\PP}$ and $S^+$.

\item{(4)} $\ooplus$ for some $\bar{\PP}$ and $S^+$ such that all elements of
$\bigcup_{\alpha<\lambda}\PP_\alpha$ are closed and, for all limit $\delta\in S^+$,
there is an $x\subseteq\delta$ with
$$\otp(x)=\cf(\delta)<\delta=\sup(x)\quad\&\quad\bwedge_{\zeta<\delta}
x\cap\zeta\in\bigcup_{\alpha<\delta}\PP_\alpha.$$

\item{(5)} $S\in I_{<f}[\lambda]$, and this is witnessed by an $S^+$ and
$\langle C_\delta:\,\delta\in S^+\rangle$ such that
for all $\delta\in S^+$:

\item{$(i)$} $\otp(C_\delta)=\cf(\delta)$ and
\item{$(ii)$} $\alpha\in\nacc(C_\delta)\implies
\alpha\in S^+\,\,\&\,\, C_\alpha=C_\delta\cap\alpha$, 

{\it Then\/}

$(1)\implies(2)\implies(3)$ and $(4)\implies(2)$ and $(5)\implies(1)$
and $(5)\implies(4).$

If we in addition assume that $f$ satisfies for all $\beta\in\lambda$,
$$\lambda\ge f(\beta)\ge\cf\bigl(f(\beta)\bigr)>\max\{\card{\beta},\aleph_0\},
\eqno{(\ast)}$$

{\it then\/} $(3)\implies(5)$ (so in this case all notions $(1)$--$(5)$
are equivalent).



{\bf Proof.} $(1)\implies(2)$ and $(5)\implies(4).$

If $S^+$ and $\langle C_\delta:\,\delta\in S^+\rangle$ witness that
$S\in I_{<f}[\lambda]$, then we can set $\bar{\PP}=\bar{\PP}[\bar{C}]$,
i.e.
$$\PP_\alpha=\{C_\beta\cap\alpha:\,\alpha\in C_\beta\},$$
and obtain $\ooplus$. We also have that all elements of $\bigcup_{\alpha<\lambda}
\PP_\alpha$ are closed.

If $S^+$
and $\bar{C}$ were chosen to witness $(5)$, then $\bar{\PP}$ witnesses $(4)$.

$(2)\implies(3)$ and $(5)\implies(1)$, as
well as $(4)\implies(2)$ are obvious.

The difficult step is to prove $(3)\implies(5)$, assuming $(\ast)$.
Luckily, the proof is exactly like the ``only if'' part of [Sh 420, 1.2].
\eop${}_{3.5.}$

\medskip

Let us from now on always assume that $f$ is a function from $\lambda$
into the cardinals, and $f(\alpha)\ge 2$ unless
stated otherwise.

We now show that $I[\lambda]$,
$I^\square[\lambda]$ and $I^{{\rm w}\square}[\lambda]$ all occur
as particular cases of $I_{<f}[\lambda]$, with the appropriate choice
of parameter $f$. Another choice of $f$ will yield the ``silly square''
of [Sh 355].

\medskip

{\bf Theorem 3.6} (1) If
$$f(\alpha)=\cases{\lambda   &if $\alpha$ is a limit ordinal\cr
                    1        & otherwise,\cr}$$
{\it then\/}
$$I_f[\lambda]=I[\lambda].$$

(2) 
$$I_1[\lambda]=I^{\square}[\lambda].$$

(3) If $f$ is given by
$$f(i)=\card{i},$$
{\it then\/}
$$I_f[\lambda]=I^{{\rm w}\square}[\lambda].$$

(4) $\lambda\setminus S^{\rm in}_\lambda\in I_\lambda[\lambda]$.

{\bf Proof.} (1) Suppose that $S\in I[\lambda]$, as witnessed by $S^+$
and $\bar{C}=\langle C_\delta:\,\delta\in \lambda\rangle.$
Then the same $S^+$, and $\langle C_\delta:\,\delta\in S^+\rangle$
witness that $S\in I_f[\lambda]$.
So, $I[\lambda]\subseteq I_f[\lambda]$.

For the other direction,
suppose that $S\in I_{f}[\lambda]$, as is
witnessed by $S^+$ and $\langle C_\delta:\,\delta\in S^+\rangle.$
We can assume that $S^+$ consists of limit ordinals only.
For $\beta\in\lambda\setminus S^+$, we define
$$C_\beta=\cases{C_\delta\cap\beta & if there is a $\delta\in S^+$ such that $
                                    \beta\in C_\delta$, and $
                                     \beta$ is a successor,\cr
                 \emptyset         & otherwise.\cr}$$
This definition is well posed, since $f(\beta)=1$ for $\beta$ a
successor.

Now $S^+$ and $\langle C_\delta:\,\delta\in\lambda\rangle$ witness
that $S\in I[\lambda]$.
We did not necessarily obtain a sequence such that ${\rm otp}(C_\delta)={\rm cf}
(\delta)$ for a club of $\delta$, but merely ${\rm otp}(C_\delta)<\delta$. It iswell known that this suffices, see [Sh 420] or [Sh 108], [Sh 88A].

(2)--(3) Both easily follow from the corresponding definitions.

(4) This is from
[Sh 355]. We simply choose for every $\alpha<\lambda$
which is a regular cardinal, a club $C_\alpha$
of $\alpha$ with
$\otp(C_\alpha)=\cf(\alpha)<\alpha$ such that
$\nacc(C_\alpha)$ contains only successor ordinals.
Of course, note that $\{\alpha<\lambda:\,\alpha=\beta^+\}$ is not stationary
in $\lambda$.\eop${}_{3.6.}$

\medskip

Just from the definition of $I_{<f}[\lambda]$, if $S\in I_{<f}[\lambda]$
is a costationary set, it seems possible that every $S^+$ which exemplifies
this satisfies that $S^+\setminus S$ is stationary. The following theorem shows that this is impossible.

\medskip

{\bf Theorem 3.7.} Suppose that
$S \in I_{<f}[\lambda]$. Then there is a set $S^+$ witnessing this, such that
$S^+= S$ modulo a non-stationary
set.

{\bf Proof.} Let us start with a set $T$ such that $T$ and $\bar\DD=
\langle D_\alpha:\,\alpha\in T\rangle$ witness that
$S\in I_{<f}[\lambda]$. Let $E$ be a club
of $\lambda$ such that $S\cap E
\subseteq T$, and we enumerate $E$ as in increasing continuous sequence
$\langle \alpha_\zeta:\,\zeta<\lambda\rangle$.

Now we set $S_1\deq\{\zeta:\,\alpha_\zeta\in E\cap S\}$,
and for $\zeta\in S_1$ we define
$C^1_\zeta\deq\{\xi<\zeta:\,\alpha_\xi\in D_{\alpha_\zeta}\}.$

Noting that the set
$$D\deq\{\zeta:\,\zeta=\alpha_\zeta\}$$
is a club of $\lambda$, we can check that $S_1=S
\,(\mod NS[\lambda])$. It is also easily seen that
$\langle C^1_\zeta:\,\zeta\in S_1\rangle$ satisfies 3.2(0) $(a)$--$(e)$, except that
it is possible that some elements of $\nacc(C^1_\zeta)$ are limit ordinals. So,
we shall define for $\zeta\in S_1$
$$C_\zeta\deq\{\epsilon:\,\epsilon\in\acc(C^1_\zeta)\}\cup
                  \{\epsilon+1:\,\epsilon\in\nacc(C^1_\zeta)\},\hbox{  if
                    }\zeta=\sup(C^1_\zeta)$$
and
$$C_{\zeta+1}\deq\{\epsilon:\,\epsilon\in\acc(C^1_\zeta)\}\cup
                  \{\epsilon+1:\,\epsilon\in\nacc(C^1_\zeta)\},\hbox{  if
                    }\zeta\neq\sup(C^1_\zeta).$$
We set $S^+\deq\{\zeta\in S^1:\,\zeta=\sup(C^1_\zeta)\}\cup
\{\zeta+1:\,\zeta>\sup(C_\zeta^1)\}$. Then $S^+$ and
$\bar C\deq\langle C_\zeta:\,\zeta\in S^+\rangle$ are as required.\eop${}_{3.7.}$

\medskip

We give a sufficient
condition
for $I_{<f}[\lambda]$ to be a normal
ideal:

\medskip

{\bf Theorem 3.8.} Suppose that $f$ satisfies for each limit $\beta\in\lambda$
$$\cf\bigl[(f(\beta))\bigr]>\beta.$$
{\it Then\/} $I_{<f}[\lambda]$ is a normal ideal on $\lambda$.

{\bf Proof.}
It can be easily checked that $I_{<f}[\lambda]$ is an ideal, we shall show that
the ideal is normal.

Suppose that $S\subseteq\lambda$ is given such that $S\notin I_{<f}[\lambda]$,
and $g$ is a regressive function on $S$. We shall assume that for all $\alpha
\in {\rm Ran} (g)$ the set $S_\alpha\deq g^{-1}(\{\alpha\})$ is an element of
$I_{<f}[\lambda]$ and obtain a contradiction.
Note that these assumptions imply that $S\cap S^{\rm in}_\lambda$ is not
stationary, so without loss of generality $S$ does not contain
any regular cardinals.

If $f(\alpha)>\lambda$ on a stationary subset $T$ of $S$, then as in the proof
of 3.6(4), we can show that $T\in I_{<f}[\lambda]$. So without loss
of generality 
(using that $I_{<f}[\lambda]$
is an ideal), $(\ast)$ of 3.5 holds on $S$. We shall show that $(3)$ of 3.5
holds for $S$, which is a contradiction.
For $\alpha\in {\rm Ran} (g)$, let $S_\alpha^+$ and $\bar{\PP^\alpha}=
\langle\PP^\alpha_\delta:\,\delta<\lambda\rangle$ be such that
${}^{-}\oplus^{f}_{{\bar\PP}^\alpha, S_\alpha^+} (S_\alpha)$ holds. For $\delta<\lambda$ we
define
$$\bar\PP_\delta\deq\bigcup_{\alpha<\delta}\{x\cup\alpha:\,x\in\PP^\alpha_\delta
\}\cup\{\alpha\}$$
and let $S^+\deq\{\delta:\,\hbox{for some }\alpha<\delta \hbox{ we have }
\delta\in S^+_\alpha\}$ and $\bar\PP\deq\langle\PP_\delta:\,\delta<\lambda
\rangle$. Then $\bar\PP$ and $S^+$ exemplify ${}^{-}\oplus_{\bar\PP,S^+}^f
(S).$\eop${}_{3.8.}$

\medskip

The family $I_{<f}[\lambda]$ does not uniquely determine $f$,
as follows from the following Fact.

\medskip

{\bf Fact 3.9.} Suppose that $f$ and $g$ are functions from
$\lambda$ into the cardinals $\ge 2$, and related to each other by the
following:
$$g(i)=\cases{ 2                  & if $i$ is not a limit ordinal, and $
                                    f(i)\le\lambda$,\cr
              \min\{f(i),\lambda^+\}  & otherwise.\cr}$$
{\it Then\/}
$$I_{<f}[\lambda]=I_{<g}[\lambda].$$

{\bf Proof.} Since $g\le f$, certainly $I_{<g}[\lambda]\subseteq
I_{<f}[\lambda].$

The other direction is simply the proof of [Sh 420, 1.2] phrased
into the language of this section.\eop${}_{3.9.}$

\medskip

\relax From the previous discussion we see that for many functions $f$, the
family $I_f[\lambda]$ is either one of the known families, or
simply obtainable as a Boolean combination of these. Now we give
a  very easy example of
a function $f$ for which this cannot be said, provided $\lambda$
is strongly inaccessible.

\medskip

{\bf Fact 3.10.} Suppose that $\lambda$ is a strongly inaccessible cardinal
and $f$ is given on $\lambda$ by
$$f(\alpha)={2^{\card{\alpha}}}^{+}.$$
{\it Then\/} $I_{<f}[\lambda]$ is a normal ideal
on $\lambda$, and $SING\cap\lambda\in I_{<f}[\lambda]$.

{\bf Proof.} That $I_{<f}[\lambda]$ is a normal ideal on $\lambda$ follows from
 3.7. 
We can choose for singular $\delta<\lambda$ any closed subset $C_\delta$
of $\delta$ with $\nacc (C_\delta)$ containing only successor ordinals,
making sure that $C_\delta$ is unbounded in $\delta$ if $\delta$ is a limit
ordinal. Then obviously for $\beta <\lambda$
$$\card{ \{C_\delta\cap\beta:\,\beta\in C_\delta\} }\le 2^{\card{\beta}}<f(\beta).$$
\centerline{
\eop${}_{3.10.}$}

\medskip

{\bf Concluding remarks 3.11.} A similar discussion could be made
with a different requirement on $\langle C_\delta:\,\delta\in S^+\rangle$
from 4.2.
We could for example require less closure and more coherence. In this
way we can recover the argument of [Sh 186, \S3], for example.
By a definition of a similar nature we can define an ideal equivalent to
$I[\lambda,\kappa)$ of [DjSh 545].

\bigskip
\eject

{\bf \S4. Appendix: More on $I[\lambda]$.} Shortly after we submitted the
(rest of the) paper for publication, we were able to prove an additional
two theorems on $I[\lambda]$, which both seem to fit with the third
section of the paper. These two theorems are the content of this appendix.

\medskip

{\bf Theorem 4.0.} Suppose that $\lambda$ is a regular cardinal,
$\kappa<\lambda$ is
regular, and for each cardinal $\sigma \in(\kappa,\lambda)$ we have:

There is a $\PP_\sigma\subseteq [\sigma]^{< \kappa}$
such that $\card{\PP_\sigma}<\lambda$ and
$$ a\in [\sigma]^{ \kappa}\implies (\exists b \in [a]^\kappa)\,([b]^{<\kappa}\subseteq\PP_\sigma).
\eqno{(\ast)}$$
{\it Then \/} there is an $S\in I[\lambda]$ such that

$$
\eqalign{
A(S)\deq\{\delta\in \lambda:\,&\cf(\delta)=\theta\hbox{ for some regular }
\theta\in(\kappa,\lambda)
\hbox{ such that}\cr 
&2^{<\kappa}<\theta\,\,\&\,\,S^\theta_\kappa\in I[\theta]
\,\,\&\,\,
S^\delta_\kappa\setminus S \hbox{ is
stationary in } \delta\}\cr}$$
is nonstationary in $\lambda$.

{\bf Proof.} The theorem will follow from Lemma 4.1.
The proof of the lemma is along the lines of various proofs presented in
[Sh 108] or [Sh 88a].

\medskip

{\bf Lemma 4.1.} Suppose that $\kappa< \lambda$ are regular cardinals
satisfying the assumptions of 4.0.
{\it Then\/} for any regular $\theta\in(\kappa,\lambda)$
such that $2^{<\kappa}<\theta$, and
for every $S_\theta\in I[\theta]$
which is a subset of $S^\theta_\kappa$, there
is an $S_\lambda=S_\lambda(\theta, S_\theta)\in I[\lambda]$ such that
$$\eqalign{B(S_\theta, S_\lambda)\deq
\bigl\{\delta\in S^\lambda_\theta:&\,\hbox{For every increasing continuous }
\langle \alpha_{\delta,i}:\,
i<\theta\rangle
\cr
&\hbox{with limit } \delta,
 \hbox { the set }\cr
&\{i<\theta:\,i\in S_\theta\hbox{ but } \alpha_{\delta,i}\notin S_\lambda\}
\hbox{ is stationary in }\theta\bigr\}\cr}$$
is not stationary in $\lambda$.

{\bf Proof of 4.1.}
Since $S_\theta\in I[\theta]$, there is a
sequence $\langle D_i:i<\theta\rangle$ and a set $S^+$ which witnesses
this, according to Definition 3.0(2). Note that the conclusion of the lemma
does not change if a nonstationary set is removed from $S_\theta$,
so we shall for convenience assume that $S_\theta\subseteq S^+$.

Let $\chi$ be large enough compared to $\lambda$, say
$\chi=\beth_9(\lambda)^+$.  
We start with an increasing continuous sequence
$\bar N=\langle N_i:\,i<\lambda\rangle$  of elementary submodels of $\langle H(\chi),\in,
\order^\ast\rangle$, which have the following properties:

\item{(a)}$\card{N_i}<\lambda.$
\item{(b)}$\{\lambda,\theta,\kappa, S_\theta,\langle D_i:\,i<\theta\rangle,
S^+,
2^{<\kappa}\}\in N_0$.
\item{(c)}$\bar N\rest (i+1)\in N_{i+1}$, for all $i<\lambda$.
\item{(d)} $N_i\cap\lambda$ is an ordinal, for all $i<\lambda$.
\item{(e)} $\theta\subseteq N_0$.

Note that $(\ast)$ as we have it, is in fact equivalent to the same statement in
which $\sigma$ is allowed to be any ordinal $<\lambda$. We shall assume this
version of $(\ast)$ for notational convenience. Then
we can without loss of generality require
that

\item{(f)} $\PP_{N_i\cap\lambda}\cup\{\PP_{N_i\cap\lambda}\}\subseteq
N_{i+1}$ for all $i<\lambda$.

Let $E\deq\{\delta<\lambda:\,N_\delta\cap\lambda=\delta\}$,
so $E$ 
consists of limit ordinals and is a club of $\lambda$ with $E\cap(\theta+1)=\emptyset$.
Define
$$\eqalign{S_\lambda\deq\{\delta\in E:\,&\cf(\delta)<\delta
\hbox{ and there is an }
A\subseteq\delta=\sup(A)\hbox{ such that}\cr
&\bigwedge_{\alpha<\delta} A\cap\alpha
\in N_\delta \,\,\&\,\,\otp(A)=\cf(\delta)\}.\cr}$$

{\it Observation.} $S_\lambda\in I[\lambda]$.

[Why?
A
more general proof is in fact given in [Sh 108], but here is a proof
using
Lemma 4.5 below. We simply
set
$\bar \RR\deq
\langle \RR_\alpha:\,\alpha<\lambda\rangle$,
where $\RR_\alpha\deq N_\alpha\cap\PP(\alpha)$.
The $\RR_\alpha$ here stand in place of $\PP_\alpha$ in Lemma 4.5.
]

{\it Proof of 4.1. continued.\/}
Assume now that
$\delta\in\acc(E)\cap S^\lambda_ \theta$.
We want to show that $\delta\notin B(S_\theta, S_\lambda)$. Let
$ \langle \alpha_{\delta,i}:\,
i<\theta\rangle$ be an increasing continuous
enumeration of $E\cap\delta$.
We shall show that $\{i<\theta:\,i\in S_\theta\hbox{ but }\alpha_{\delta,i}
\notin S_\lambda\}$ is not stationary in $\theta$.

\smallskip

{\it Observation.\/}
Since $S_\theta\in I[\theta]$,
we can find a sequence $\langle C_{\alpha_{\delta,i}}:\,i<\theta\rangle$ such that:

\item{-} $C_{\alpha_{\delta,i}}$ is a subset of $\alpha_{\delta,i}$.
\item{-} $\otp(C_{\alpha_{\delta,i}})\le \kappa.$
\item{-} $\beta\in C_{\alpha_{\delta,i}}\implies
\beta=\alpha_{\delta,j}\hbox{ for some } j$ and $C_\beta
=C_{\alpha_{\delta,i}}\cap \alpha_{\delta,j}.$
\item{-} $i\in S_\theta\implies \alpha_{\delta,i}=\sup (C_{\alpha_{\delta,i}}).$

[How do we find such a sequence? First
we set $C^\ast_i$ to be the first $\kappa$ nonaccumulation
points of $D_i$,
for $i<\theta$.
Then let $C_{\alpha_{\delta,i}}\deq\{\alpha_{\delta,j}:\,j\in C^\ast_i\}.$]

{\it Proof of 4.1 continued.\/} Note that the sequence $\langle C_i^\ast:\,i<\theta
\rangle$ is both an element and a subset of $N_0$.
We also have that every $C_{\alpha_{\delta,i}}$ is in $N_\delta$,
but note that we do not know that necessarily $C_{\alpha_{\delta,i}}
\in N_{\alpha_{\delta,i}+1}.$ So we shall define a
function $h:\theta\to \theta$
by
$$h(i)\deq\min\bigl\{\epsilon<\theta:
\,(\forall x\in [C_i^\ast]^{<\kappa})\,(\{\alpha_{\delta,j}:\,j\in x\}
\in N_\delta\implies \{\alpha_{\delta,j}:\,j\in x\}\in
 N_{\alpha_{\delta,\epsilon}})\bigr\}.$$
Note that if $x\in [C_i^\ast]^{<\kappa}$
and $\{\alpha_{\delta,j}:\,j\in x\}\in
N_\delta$, then $\{\alpha_{\delta,j}:\,j\in x\}\in N_{\alpha_{\delta,\epsilon}}$
for some $\epsilon<\theta$,
as $2^{<\kappa}<\theta$. (So  $h(i)$ is well defined for all $i<\theta$.)

\noindent{[Why? We can find a $\xi<\delta$
such that $N_\xi$ already contains all bounded subsets of $\kappa$
that are going to appear in $N_\delta$, as otherwise we would be able to
inductively construct $\theta$ many bounded subsets of $\kappa$, in
contradiction with $2^{<\kappa}<\theta$.]}

Let
$$e\deq\{i<\theta:\,(\forall j<i)
(h(j)<i)\}.$$
Obviously, $e$ is a club of $\theta$.
We claim $S_\theta\cap e\subseteq\{i<\theta:\,\alpha_{\delta,i}\in
S_\lambda\}$. (Hence, $\{i<\theta:\,i\in S_\theta\hbox{ but }
\alpha_{\delta,i}\notin S_\lambda\}$ is not stationary in $\lambda$.)

To see this,
consider an $\alpha_{\delta,i}$ for some
$i\in S_\theta\cap e$.
We know that $\alpha_{\delta,i}\in E$ and $\cf (\alpha_{\delta,i})
=\cf (i)=\kappa <\alpha_{\delta,i}$, as $i\in S_\theta$.
Now, $C_{\alpha_{\delta,i}}$ is unbounded in 
$\alpha_{\delta,i}$, so $\otp(C_{\alpha_{\delta,i}})=\kappa$. Therefore
for some $A_i\subseteq C_{\alpha_{\delta,i}}$
we  have 
$$\otp(A_i)=\kappa\,\,\&\,\,[A_i]^{<\kappa}\subseteq \PP_{\alpha_{\delta,i}}.$$
So
$\sup(A_i)=\alpha_{\delta,i}$, and it suffices
to see that for all $\alpha<\alpha_{\delta,i}$ we have $A_i\cap\alpha
\in N_{\alpha_{\delta,i}}$.

Let us fix an $\alpha<\alpha_{\delta,i}$. Let $j^\ast=\min\{j<i:\,
\alpha\le\alpha_{\delta,j}\in A_i\}.$ Therefore
$$A_i\cap \alpha
=A_i\cap \alpha_{\delta, j^\ast}=A_i\cap C_{\alpha_{\delta,i}}\cap
\alpha_{\delta,j^\ast}=A_i\cap C_{\alpha_{\delta,j^\ast}}.$$
Let $x\subseteq C_{j^\ast}^\ast$ be such that
$A_i\cap C_{\alpha_{\delta,j^\ast}}=\{\alpha_{\delta,j}:\, j\in x\}, $ so
$\card{x}<\kappa.$
Note that by the choice of $A_i$ we have
$A_i\cap \alpha\in \PP_{\alpha_{\delta,i}}.$
Since $\PP_{\alpha_{\delta,i}}\subseteq N_{\alpha_{\delta,i}+1}\subseteq
N_\delta$, we have $A_i\cap \alpha=\{\alpha_{\delta,j}:\, j\in x\}\in
N_\delta$, therefore $A_i\cap \alpha\in N_{\alpha_\delta, h(j)}$.
As $h(j)<i$,
our claim is correct.

Consequently, $\delta\notin B(S_\theta, S_\lambda)$, so $B(S_\theta,
S_\lambda)\cap \acc(E)=\emptyset$,
hence $B(S_\theta, S_\lambda)$
is not stationary.
\eop${}_{4.1.}$

{\bf Proof of 4.0. continued.} Suppose that
$\kappa$ and $\lambda$ are as in the assumptions of the theorem, and
$\theta\in (\kappa,\lambda)$ is regular and such that
$2^{<\kappa}<\theta$ and $S^\theta_\kappa\in I[\theta]$.
We will apply Lemma 4.1 to $\kappa,\theta$ and $\lambda$.
By the Lemma, we can find a set $S_\lambda(\theta)\in
I[\lambda]$ such that $B(S^\theta_\kappa, S_\lambda(\theta))$ is
not stationary in $\lambda$.

Let $S$ be the diagonal union
of $\{S_\lambda(\theta):\, S_\lambda(\theta)\hbox{ is defined}\}$,
i.e.
$$S\deq\bigl\{\alpha<\lambda:\,\bigl(\exists\theta<\alpha
\bigr)\bigl(S_\lambda(\theta)
\hbox{ is defined and }\alpha\in S_\lambda(\theta)\bigr)\bigr\}.$$
It follows from the normality 
of $I[\lambda]$ that it is closed
under formation of such unions, so $S\in I[\lambda]$. 

Now suppose that $A(S)$ is stationary in $\lambda$, then there must be
a regular $\theta<\lambda$
such that $2^{<\kappa}<\theta$
and $S^\theta_\kappa\in I[\theta]$, and
$$\{\delta\in S^\lambda_\theta:\,S^\delta_\kappa\setminus S
\hbox{ is stationary in }\delta\}$$
is stationary in $\lambda$.
Note that $S_\lambda(\theta)$ is defined. Now,
for any $\delta>\theta$ in the above set, the set $S^\delta_\kappa
\setminus S_\lambda(\theta)$ is stationary
in $\delta$. In particular for every sequence $\langle \alpha_{\delta,i}:
\,i<\theta\rangle$ which increasingly enumerates a club of $\delta$,
the set 
$$\{i<\theta:\,\cf(i)=\kappa\,\,\&\,\,\alpha_{\delta, i}\notin S_\lambda(
\theta)\}$$
is stationary in $\theta$.
Hence $\delta\in B(S^\theta_\kappa,
S_\lambda(\theta))$, and the set of such $\delta$ is nonstationary
in $\lambda$, by the Lemma.
This is a contradiction, hence $S$ is as required.\eop${}_{4.0.}$

\medskip

{\bf Remark and Conclusion 4.2.} Property $(\ast)$ was considered in [Sh 430].
Obviously, $(\ast)$ is implied
by
$$(\forall\sigma<\lambda)\,(\sigma^{<\kappa}<\lambda).$$
Also, $S^\theta_\kappa\in I[\theta]$ obviously follows
from $(\forall\alpha<\theta)(\card{\alpha}^\kappa<\theta).$
Hence, if for example $\lambda\ge\beth_\omega$ and $\kappa<\beth_\omega$
are regular, there is a set $S\in I[\lambda]$ and an
$n\in \omega$ such that for all
regular $\theta\in (\beth_n,\lambda)$
$$\{\delta\in S^\lambda_\theta:\,S^\delta_\kappa\setminus S
\hbox{ is stationary in }\delta\}$$
is nonstationary in $\lambda$.



\medskip   

{\bf Definition 4.3.} Let $\lambda$  be a regular uncountable cardinal.

 (1) For $\Theta\subseteq REG$, let
$$S_{\Theta}\deq\{\delta:\,\delta>\cf(\delta)\in \Theta\}.$$

(2) A sequence $\langle {\CC}_\alpha:\,\alpha <\lambda\rangle$ is said to be
a $\Theta$-{\it weak square on\/} $S\subseteq \lambda$ iff
there is an $S^+\subseteq\lambda$ such that
$S\setminus S^+$ is nonstationary, and for every $\alpha<\lambda$:

\item{$(i)$} $\CC_\alpha$ is a nonempty family of subsets of $\alpha$.
\item{$(ii)$} $\card{\CC_\alpha} <\lambda.$
\item{$(iii)$} $C\in \CC_\alpha \implies C$ is closed.
\item{$(iv)$} $\emptyset\neq
C\in \CC_\alpha\,\,\&\,\,\alpha\hbox{ is a limit ordinal}\implies
\alpha\in S^+\,\,\&\,\,
\sup(C)=\alpha\,\,\&\,\,\otp(C)<\min(C).$
\item{$(v)$} $\beta\in C\,\,\&\,\, C\in\CC_\alpha\,\,
\&\,\,\beta\in S_\Theta
\implies \beta\in S^+\,\,\&\,\, C\cap\beta\in \CC_\beta.$

\medskip

{\bf Theorem 4.4.} Let $\lambda$  be a regular uncountable cardinal.
Suppose that $S\in I[\lambda]$ and
$$\Theta\deq\{\theta<\lambda:\,\cf(\theta)=\theta\,\,\&\,\,\hbox{every tree with }
<\lambda \hbox{ nodes has }<\lambda \hbox{ branches of length }\theta\}.$$
{\it Then\/} there is a $\Theta$-weak square on $S$.

\medskip

The following Lemma was proved in  [Sh 420]:

\medskip

{\bf Lemma 4.5.}  Let us define $I'[\lambda]$ as the family of all subsets $S$ of $\lambda$
for which there is an $S^+$ and $\langle \PP_\alpha:\,\alpha<\lambda\rangle$
such that $S\setminus S^+$ is nonstationary and
for every $\alpha\in\lambda$:

\item{(A)} $\PP_\alpha$ is a family of $<\lambda$ subsets of $\alpha$.
\item{(B)} If $\alpha\in S^+$, then there is an unbounded
subset $a\subseteq \alpha$ such that
$$(\forall\beta\in\alpha)(a\cap\beta\in \PP_\beta).$$
{\it Then\/}:

\item{(1)} $I'[\lambda]=I[\lambda].$

\item{(2)} Without loss of generality we can require in the definition of $I^{'}[\lambda]$:

\item{(C)} The sets in $\cup_{\alpha<\lambda}\PP_\alpha$ are closed
and for $\alpha\in S^+$, the set $a$ from (B) satisfies
$$\otp(a)=\cf(a)<\min(a).$$

\smallskip
 
We shall for completeness include the proof of the first part of this lemma.
We do not prove (2), as the proof is the same as that of 3.5.

{\bf Proof of 4.5.}(1) It follows by 3.5 and 3.6(1) that 
$I'[\lambda]\subseteq I[\lambda]$, so let us start with an
$S\in I[\lambda]$ and show that $S\in I'[\lambda]$.
Let $\bar\PP_\alpha=\langle\PP^0_\alpha:\,\alpha<\lambda\rangle$ and $S^+$ exemplify that
$S\in I[\lambda]$.(We use the definition from [Sh 108]).

For $\alpha<\lambda$ let us define $\PP^1_\alpha$ by
$$\PP^1_\alpha\deq\{ c\cap[\gamma,\beta):\,\gamma\le\beta\le\alpha\,\,\&\,\,
c\in \PP^0_\alpha\},$$
and let
$$\PP_\alpha^2\deq
\bigcup_{\beta\le\alpha}\PP_\beta^1.$$
For every $a\in \bigcup_{\alpha<\lambda}\PP_\alpha^2$, we define a function
$f_a$ with ${\rm Dom}(f_a)=a$ by
$$f_a(\gamma)\deq\min\{\beta\ge\gamma:\,a\cap\gamma\in \PP_\beta^2\}.$$
Note that $f_a(\gamma)$ is well defined for any $\gamma\in 
a$.
Finally, let
$$\PP_\alpha\deq\{{\rm{Ran}}(f_a)\cap [\gamma,\beta):\,a\in\PP_\alpha^2\,\,
\&\,\,\gamma\le\beta\le\alpha\}.$$

Let us check that $\bar\PP\deq\langle\PP_\alpha:\,\alpha<\lambda\rangle$ and $S^+$
exemplify that $S$ is in $I'[\lambda]$.

First we need to observe the following facts:

\item{(1)} $\gamma\le f_a(\gamma)$, if $\gamma\in a \in \bigcup_{\alpha<\lambda}\PP_\lambda^2.$

This follows just from the definition of $f_a$.

\item{(2)} $\gamma_1\le\gamma_2\in a \in \bigcup_{\alpha<\lambda}\PP_\lambda^2\implies
f_a(\gamma_1)\le f_a(\gamma_2).$

Let $\beta_1\deq f_a(\gamma_1)$. Suppose that $f_a(\gamma_2)<\beta_1$, so in particular
there is a $\beta_2<\beta_1$ such that $a\cap\gamma_2\in\PP_{\beta_2}^2$. By the
definition of $\PP_{\beta_2}^2$, there is some
$\beta_3\le\beta_2$ such that
$a\cap\gamma_2\in\PP_{\beta_3}^1.$ Therefore $a\cap\gamma_2=b\cap [\beta_4,\beta_5)$
for some $\beta_4\le\beta_5\le\beta_3$ and $b\in \PP_{\beta_3}^0$.
In particular $\gamma_1\in b\cap [\beta_4,\beta_5)$. Hence $b\cap [\beta_4,
\gamma_1)=a\cap\gamma_1$ is in $ \PP_{\beta_3}^1$, so in $\PP_{\beta_2}^2$,
which contradicts the minimality of $\beta_1$ in the definition of $f_a(\gamma_1)$.

\item{(3)} $\gamma\in a\implies f_{a\cap\gamma}=f_a\rest(a\cap\gamma)$.

First note that $f_{a\cap\gamma}$ is well defined an then just use the definition of $f_a$.

Now we check that $\bar\PP$ has the required properties. It is obvious from
the definition of $\bar\PP$ that
$\alpha<\lambda\implies\PP_\alpha\subseteq\PP(\alpha)$ and
$\card{\PP_\alpha}<\lambda$.
Now suppose that $\alpha\in S^+$, and let $a$ be an unbounded subset of $\alpha$
with $\otp(a)<\alpha$ and
$$(\forall \beta<\alpha)(a\cap\beta\in \bigcup_{\gamma<\alpha}\PP_\gamma^0).$$
Notice that for every $\gamma<\lambda$ we have $\PP_\gamma^0\subseteq
\PP_\gamma^1\subseteq\PP_\gamma^2$, so $f_{a}$ is defined.
By the choice of $a$ and $\alpha$ we see that ${\rm Ran}(f_{a})
\subseteq \alpha$. Moreover, by (1) above, ${\rm Ran}(f_{a})$ is unbounded
in $\alpha$. Since ${\rm Dom}(f_{a})=a$, then $\otp({\rm Ran}(f_a))\le\otp
({\rm Dom})(f_a)<\alpha$.

It remains to see that if $\beta\in {\rm Ran}(f_a)$, then $\beta\cap {\rm Ran}(f_a)
\in\PP_\beta$. By (1)+(3) above we know that ${\rm Ran}(f_a)\cap\beta
={\rm Ran}(f_{a\cap\beta})$, so we have proved our Lemma.\eop${}_{4.5(1)}$

\medskip
  
{\bf Proof of Theorem 4.4.}
 Let $S^+$ and $\langle\PP_\alpha:\,\alpha<\lambda\rangle$ exemplify
that $S\in I[\lambda]$. By the above Lemma we can assume that for every $\alpha\in S^+$,
there is an unbounded subset $a$ of $\alpha$ such that
$$(\forall\beta<\alpha)\,(a\cap\beta\in\PP_\beta)
\,\,\&\,\,\otp(a)=\cf(a)<\min(a)$$
and that each $\PP_\alpha$ consists of closed sets.
Let us define
$$\CC_\alpha=\cases{\{a\subseteq\alpha:\,
(\forall\beta<\alpha)\,(a\cap\beta\in\PP_\beta)\cr
\,\,\&\,\,a\neq\emptyset\implies
\cf(\alpha)=\otp(a)<\min(a)\cr
\,\,\&\,\,
\alpha\hbox{ limit}\implies\sup(a)=\alpha\}\cr
 &if $\alpha\in S_\Theta
\cap S^+$ or\cr
&$\alpha$ a successor\cr
\{\emptyset\} &otherwise.\cr}$$

Then it is easy to check all the requirements for a $\Theta$-weak square.
That $\CC_\alpha$ is never empty follows from the definition of $I[\lambda]$.
To see that for $\alpha\in S_\Theta$ we have $\card{\CC_\alpha}<\lambda$,
consider the tree $T$ which is defined in the following way.

For $\beta<\alpha$, we have that the $\beta$-level of $T$ is
$${\rm lev}_\beta(T)\deq\{a\cap\beta:\,a\in \CC_\alpha\}$$
and $T$ is ordered by $\subseteq$. Then $T$ has $<\lambda$ nodes,
and every element of $\CC_\alpha$ is the union of an $\alpha$-branch
of $T$.
As $\theta\deq\cf(\alpha)\in \Theta$, 
by the definition of $\Theta$ we have that
$\card{\CC_\alpha}<\lambda$.

Also notice that all elements of $\CC_\alpha$ are increasing unions
of closed sets without the last element, so they are closed.

Finally, if $\beta\in C\,\,\&\,\, C\in\CC_\alpha
\,\,\&\,\,\beta\in S_\Theta$, then
$\beta\in S^+$ and $C\cap\beta\in \CC_\beta.$
\eop${}_{4.4.}$

\medskip

{\bf Remark 4.6.} If $\lambda>\mu$
is regular and $\mu$ is a strong limit singular with
$\kappa=\cf(\mu)$, then $\Theta=REG\cap\mu\setminus\{\kappa\}$ satisfies
the condition of Theorem 4.4, by [Sh 460 1.1].
 
\eject
\baselineskip=12pt

\bigskip
\bigskip

\centerline{\bf REFERENCES}

\bigskip
\bigskip
\bigskip

\item{[DjSh 545]} Mirna D\v zamonja and Saharon Shelah, Saturated filters at
successors of singulars, weak reflection and 
yet another weak club principle, {\it submitted.}

\item{[Gi]} Moti Gitik, Some results on nonstationary ideal,
{\it Israel J. of Math., to appear}

\item{[Gi1]} Moti Gitik, Some results on nonstationary ideal II, {\it 
Israel J. of Math.,submitted.}

\item{[Sh -b]} Saharon Shelah, Proper Forcing,
{\it Lect. Notes in Mathematics\/} 940 (1982), {\it Springer-Verlag.}

\item{[Sh -e]} Saharon Shelah,
Non structure theory, {\it accepted, Oxford University Press.}

\item{[Sh -f]} Saharon Shelah, Proper and improper forcing,
{\it accepted, Oxford University Press.}

\item{[Sh -g]} Saharon Shelah, Cardinal Arithmetic, 
{\it Oxford University Press\/} 1994.

\item{[Sh 88a]} Appendix: on stationary sets (to ``Classification of
nonelementary classes II. Abstract elementary classes''), in
{\it Classification theory (Chicago, IL, 1985), Proceedings
of the USA-Israel Conference on Classification Theory, Chicago, December
1985; J.~T. Baldwin (ed.)\/}, Lecture Notes in Mathematics, (1292)
483--495, {\it Springer-Verlag\/} 1987.

\item{[Sh 108]} Saharon Shelah, On Successors of Singular Cardinals, in
{\it Logic Colloquium 78, M. Boffa, D. van Dalen, K. McAloon (eds.)\/},
357-380,
{\it North-Holland Publishing Company\/} 1979.

\item{[Sh 186]}
Saharon Shelah, Diamonds, Uniformization, {\it Journal of Symbolic
Logic \/} (49) 1022-1033, 1984.

\item{[Sh 237 e]} Remarks on Squares, in {\it Around Classification Theory
of Models\/}, Lecture Notes in Mathematics, (1182) 276--279,
{\it Springer-Verlag\/} 1986.

\item{[Sh 351]} Reflecting stationary sets and successors of singular
cardinals, {\it Archive for Mathematical Logic\/} (31) 1--29, 1991.

\item{[Sh 355]} Saharon Shelah, $\aleph_{\omega+1}$ has a Jonsson Algebra, in
{\it Cardinal Arithmetic, Chapter II\/}, {\it Oxford University
Press\/} 1994.

\item{[Sh 365]} Saharon Shelah, Jonsson Algebras in Inaccessible Cardinals, in
{\it Cardinal Arithmetic, Chapter III\/}, {\it Oxford University
Press\/} 1994.

\item{[Sh 400]} Saharon Shelah, Cardinal Arithmetic, in {\it Cardinal
Arithmetic, Chapter IX\/}, {\it Oxford University Press\/} 1994.

\item{[Sh 410]} Saharon Shelah, More on Cardinal arithmetic,
{\it Archive for Mathematical Logic\/} (32) 399--428, 1993.

\item{[Sh 420]} Saharon Shelah, Advances in Cardinal Arithmetic, in {\it
Proceedings of the Banff Conference in Alberta 4/91\/}, accepted.

\item{[Sh 430]} Saharon Shelah, Further Cardinal Arithmetic,
{\it Israel Journal of Mathematics\/}, accepted. 

\item{[Sh 460]} Saharon Shelah, The Generalized Continuum Hypothesis
revisited, {\it Israel Journal of Mathematics\/}, submitted.

\item{[Sh 580]} Saharon Shelah, Strong covering revisited, {\it
to appear in a special issue of Fundamenta Mathematicae.}

\eject
\baselineskip=24pt
\bye